\documentclass[10pt]{article}
\usepackage[latin1]{inputenc}
\usepackage[english]{babel}
\usepackage{amssymb,amsfonts, amsmath}
\newtheorem{theorem}{Theorem}[section]
\newtheorem{proposition}[theorem]{Proposition}

\newtheorem{lemma}[theorem]{Lemma}
\newtheorem{remark}[theorem]{Remark}

\newcommand{\bcl}{\begin{center}}
\newcommand{\ecl}{\end{center}}
\newcommand{\brl}{\begin{right}}
\newcommand{\erl}{\end{right}}
\newcommand{\ben}{\begin{enumerate}}
\newcommand{\een}{\end{enumerate}}
\newcommand{\barr}{\begin{array}}
\newcommand{\earr}{\end{array}}
\newcommand{\btab}{\begin{tabular}}
\newcommand{\etab}{\end{tabular}}
\newcommand{\bdoc}{\begin{document}}
\newcommand{\edoc}{\end{document}}
\newcommand{\beqy}{\begin{eqnarray}}
\newcommand{\eeqy}{\end{eqnarray}}

\newcommand{\beqi}{\begin{eqnarray*}}
\newcommand{\eeqi}{\end{eqnarray*}}
\newcommand{\bitem}{\begin{itemize}}

\newcommand{\eitem}{\end{itemize}}
\newcommand{\nln}{\newline}
\newcommand{\newt}{\newtheorem}
\newcommand{\pa}{\partial}
\newcommand{\re}{{I\!\!R}}
\newcommand{\ren}{\re^N}
\newcommand{\xr}{x\in\re }
\newcommand{\x}{\times}
\newcommand{\dyle}{\displaystyle}
\newcommand{\ene}{{I\!\!N}}
\newcommand{\irn}{\int\limits_{\re^N}}
\newcommand{\io}{\int\limits_{\O}}
\newcommand{\meas}{{\rm meas\,}}

\newcommand{\sign}{{\rm sign}}
\newcommand{\map}{\longrightarrow }
\newcommand{\imp}{\Longrightarrow }
\newcommand{\sen}{{\rm sen\,}}
\newcommand{\tg}{{\rm tg\,}}
\newcommand{\arcsen}{{\rm arcsen\,}}
\newcommand{\arctg}{{\rm arctg\,}}
\newcommand{\supp}{{\textsl supp\ }}
\newcommand{\ity}{\int_{-\iy}^{+\iy}}
\newcommand{\limit}{\lim\limits}
\newcommand{\limi}{\limit_{n\to\infty}}
\newcommand{\sumi}{\sum\limits_{n=1}^{\infty}}
\newcommand{\ulu}{\underline u}
\newcommand{\ulw}{\underline w}
\newcommand{\ulz}{\underline z}
\newcommand{\ulv}{\underline v}
\newcommand{\uls}{\underline s}
\newcommand{\olu}{\overline u}
\newcommand{\olv}{\overline v}
\newcommand{\ols}{\overline s}
\newcommand{\ob}{\overline\b}
\newcommand{\ovar}{\overline\var}
\newcommand{\wv}{\widetilde v}
\newcommand{\wu}{\widetilde u}
\newcommand{\ws}{\widetilde s}
\renewcommand{\a }{\alpha }
\renewcommand{\b }{\beta }
\newcommand{\g }{\gamma}
\newcommand{\G }{\Gamma }
\renewcommand{\d }{\delta }

\newcommand{\D }{\Delta }
\newcommand{\e }{\varepsilon }
\newcommand{\z }{\zeta }
\renewcommand{\l }{\lambda }
\renewcommand{\L }{\Lambda }
\newcommand{\m }{\mu }
\newcommand{\n }{\nabla }
\newcommand{\s }{\sigma }
\newcommand{\Sig }{\Sigma }
\renewcommand{\t }{\tau }
\newcommand{\var }{\varphi }
\renewcommand{\o }{\omega }
\renewcommand{\O }{\Omega }
\newcommand{\bR}{{\bf R}}
\newcommand{\bC}{{\bf C}}
\newcommand{\bZ}{{\bf Z}}
\newcommand{\bN}{{\bf N}}
\newcommand{\bQ}{{\bf Q}}
\newcommand{\bK}{{\bf K}}
\newcommand{\bI}{{\bf I}}
\newcommand{\bv}{{\bf v}}
\newcommand{\bV}{{\bf V}}
\DeclareMathOperator{\di}{div} \DeclareMathOperator{\Hess}{Hess}
\DeclareMathOperator{\conv}{conv} \DeclareMathOperator{\cd}{D}
\DeclareMathOperator{\ud}{\overline{D}}
\DeclareMathOperator{\suppo}{supp} \DeclareMathOperator{\Ric}{Ric}
\DeclareMathOperator{\del}{\delta} \DeclareMathOperator{\inj}{inj}

\def\qed{\unskip\kern 6pt \penalty 500
\raise -2pt\hbox{\vrule \vbox to10pt{\hrule width 4pt
\vfill\hrule}\vrule}\par}

\newenvironment{Proof}{\removelastskip\vskip12pt
plus 1pt \noindent\em\rm}{\hfill {\qed \hskip .2cm}}

\title
{Allen-Cahn Approximation of Mean Curvature Flow in Riemannian
manifolds I,\\
uniform estimates}

\author{Adriano Pisante\thanks{Dipartimento di Matematica "G.
Castelnuovo", Universit\`a di Roma ``La Sapienza'', P.le A.~Moro
5, I-00185 Roma, Italia (pisante@mat.uniroma1.it).} \ \ and Fabio
Punzo\thanks{Dipartimento di Matematica "G. Castelnuovo",
Universit\`a di Roma ``La Sapienza'', P.le A.~Moro 5, I-00185
Roma, Italia (punzo@mat.uniroma1.it).}}

\date{}

\begin{document}
 \maketitle
 {\abstract{ \noindent We are concerned with solutions to the
parabolic Allen-Cahn equation in Riemannian manifolds. For a general class of initial condition
we show non positivity of the limiting energy discrepancy. This in turn allows
to prove
almost monotonicity formula (a weak counterpart of Huisken's
monotonicity formula) which gives a local uniform control of the energy densities at small scales.

Such results will be used in \cite{PiPu2} to extend previous
important results from \cite{Ilm1} in Euclidean space, showing
convergence of solutions to the parabolic Allen-Cahn equations to
Brakke's motion by mean curvature in space forms.

\bigskip

%\noindent {\it  2010 Mathematics Subject Classification: .}

}}

\bigskip
\medskip
\smallskip

\section{Introduction} \setcounter{equation}{0}
We are concerned with the Allen-Cahn equation
\begin{equation}\label{e1}
\pa_t u^\e \,=\, \Delta u^\e -\frac 1{\e^2} f(u^\e)\quad
\textrm{in}\;\; M\times(0,\infty)\,,
\end{equation}
completed with the initial condition
\begin{equation}\label{e1a}
u^\e=u^\e_0\quad \textrm{in}\;\; M\times \{0\}\,.
\end{equation}
Here $\e>0$ is a small parameter, $M$ is an $N-$dimensional
Riemannian manifold with Ricci curvature bounded from below,
$\Delta$ is the Laplace-Beltrami operator on $M$, the function $f$
is the derivative of a potential $F$ with two wells of equal depth
at $u = -1$ and at $u=1$. To be specific, we will always assume
for simplicity that $f$ satisfies
\[
\textrm{\ \ } \left\{
\begin{array}{l}
(i) \quad\,\; f = F',\; \textrm{with}\;\; F\in C^\infty (\re),\,
F\;\,\textrm{even}\,;
\\
(ii) \quad f(0)=f(\pm1)=0\,, f<0\,\,\hbox{in}\,\, (0,1)\,, f>0\,\,\hbox{in}\,\, (1,\infty),\,\\ \qquad \;\, f'(0)<0,
f'(\pm1)>0\,;
\\ (iii)\;\; \,  F>0\;\,\textrm{in}\;\, \re\setminus\{\pm1\},\; F(\pm1)=0\,;
\\ (iv) \;\;\, \min_{[\a, \infty)}F''>0,\, \textrm{for some}\;\, \a\in
(0,1)\,.
\end{array}
\right.\leqno(H_0)\] A typical example is
\[F(u)=\frac 1 2 (1-u^2)^2,\,\, f(u)=2u(u^2-1)\,.\]
We set
$$f_\e(u):= \frac 1{\e^2} f(u)\,,\quad F_\e(u):=\frac 1{\e^2}F(u)\,.$$

\smallskip

Observe that problem
\begin{equation}\label{e1i}
\pa_t u^\e \,=\, \Delta u^\e -\frac 1{\e^2} f(u^\e)\quad
\textrm{in}\;\; \re^N\times(0,\infty)\,,
\end{equation}
\begin{equation}\label{e2i}
u^\e=u^\e_0\quad \textrm{in}\;\; \re^N\times \{0\},
\end{equation}
which corresponds to problem \eqref{e1}-\eqref{e1a} in the special
case $M=\re^N$, has been the object of detailed investigations in
order to describe formation and evolution of interfaces
(approximatively) driven by their mean curvature. Indeed, it is
well-known that the term $f_\e$ forces the solution $u^\e$ to
problem \eqref{e1i}-\eqref{e2i} to problem to take values $1$ or
$-1$, as $\e\to 0^+$; moreover, the interface that separates the
two regions of $\re^N\times (0,\infty)$ in which $u^\e$ converges
to $1$ or $-1$, say the region where $\{|u^\e|<\frac 1 2\}$, is a
set of thickness of order $\e$ that in the limit as $\e\to 0$
approximatively moves by mean curvature flow as long as time
varies. A large number of papers have been devoted to this type of
results, using several methods. Without any claim for
completeness, we mention \cite{AHM}, \cite{BSS},
\cite{BK1}-\cite{dMS2}, \cite{ESS}-\cite{Fu}, \cite{Ilm1},
\cite{RSK}, \cite{Son}; observe that also similar questions have
been addressed also for the stationary equation (see, $e.g.$,
\cite{HT}, \cite{Modica}, \cite{PacRit}, \cite{PisPons}) and for
systems (see, $e.g.$, \cite{BOS}, \cite{JerrSon}). For a
comprehensive account of literature on this subject, also
containing the description of main results obtained and various
methods used, we refer the reader to \cite{Son} and references
therein. In this connection, note that in the literature several
notions of mean curvature flows have been considered (see, $e.g.$,
\cite{BellNov}, \cite{BellPaol}, \cite{CGG}, \cite{ES2}-\cite{ES4}, \cite{ESS},
\cite{Ilm2}, \cite{Nov}, \cite{Son}).

In the sequel, before describing the results of the present paper,
we limit ourselves to recall those established in \cite{Ilm1};
moreover, in general, we shall briefly explain the line of
arguments followed there to obtain them.

%To this purpose, let us recall that if $\{M_t\}_{t\geq 0}$ is a family of submanifolds of $\re^N$ that evolve by mean curvature
%flow associated
%to the mean curvature vector $\vec H$, then for any test function $\phi\in C^\infty_c(\re^N), \phi\geq 0$ there holds
%\begin{equation}\label{e3i}
%\frac{d}{dt}\int_{M_t} \phi \,d\mathcal H^{N-1} \leq  \int_{M_t} \Big\{ -\phi H^2  + \langle \nabla \phi, \vec H  \rangle \Big\}\, d\mathcal H^{N-1}\,.
%\end{equation}
%This inequality is called {\it Brakke's inequality} (see
%\cite{Brakke}); moreover, usually the first term in the right hand
%side is referred to as the {\it shrinkage} term, while the second
%one as the {\it transport} term.

An important role is played by the one-dimensional standing wave
$q^\e$  for \eqref{e1i}, for which there hold:
\begin{equation}\label{e6i}
q^\e_{rr}(r)-f_\e(q^\e(r))\,=\,0\,, \quad r\in \re\,,
\end{equation}
\begin{equation}\label{e7i}
q^\e_r>0, q^\e(0)=0, \lim_{r\to +\infty}q^\e(r)=1, \lim_{r\to
-\infty}q^\e(r)=-1\,.
\end{equation}
Concerning initial conditions, as model case one usually considers
well prepared data, i.e. in the form:
\begin{equation}\label{u70}
u_0^{\e}(x):= q^\e(\tilde d(x, \Sigma_0))\, \quad
(x\in M)\,,\end{equation}
where $\tilde d(\cdot, \Sigma_0)$ is the
{\it signed distance} from the smooth boundary $\Sigma_0$ of a
bounded domain $E_0\subset \re^N$, defined by
\[\tilde d(x,\Sigma_0):= d(x, \Sigma_0)\quad \textrm{if}\;\, x\in E_0\,,\quad \tilde d(x,\Sigma_0):= - d(x,\Sigma_0)\quad \textrm{if}\;\, x\in \re^N\setminus E_0\, \, ,\]
possibly regularizing the signed distance far away from the initial interface $\Sigma_0$.

In \cite{Ilm1} it is proved that when $\e>0$ is sufficiently small
and $u^\e$ solves problem \eqref{e1i}-\eqref{e2i}, the {\it energy
density}
\[d\mu^\e_t := \left\{\frac{\e}2 |\nabla u^\e|^2 + \frac 1{\e}F(u^\e)\right\}dx,\]
satisfies both a {\it Brakke's} and a {\it Huisken's} type
formula, in analogy to Brakke's inequality and Huisken's
monotonicity formula for a family $\{\Sigma_t\}_{t\geq 0}$ of
hypersurfaces of $\re^N$ that evolve by mean curvature flow.
However, such formulas for $d\mu^\e_t$ involve a new term: the
{\it discrepancy} Radon measure
\[d\xi^\e_t:=\left\{\frac{\e}2 |\nabla u^\e|^2 - \frac 1{\e}F(u^\e)\right\}dx\,.\]

A crucial point in \cite{Ilm1} is to show that
\begin{equation}\label{e4i}
d\xi^\e_t \, \leq \,0\,\quad \textrm{for all}\;\; t\geq 0, \e>0\,.
\end{equation}
This inequality, following \cite{Modica} in the stationary case,
is deduced from the inequality
\begin{equation}\label{e5i}
|\nabla r^\e (x,t)|\leq 1\quad \textrm{for all}\;\; x\in \re^N,
t\geq 0,
\end{equation}
where $r^\e:\re^N\times (0,\infty)$ is the function defined by
\[u^\e(x,t)=q^\e(r^\e(x,t))\quad \textrm{for all}\;\; x\in \re^N, t>0.\]
Note that, by hypothesis \eqref{u70} on initial conditions
$u_0^\e$, inequality \eqref{e5i} is satisfied for all $x\in \re^N$
and $t=0$, since $x\mapsto \tilde d(x, \Sigma_0)$ is
$1-$Lipschitz. Then by maximum principle, applied to a certain
parabolic equation satisfied by $z:= |\nabla r^\e|^2$, it is
obtained for all $x\in \re^N, t>0$.

As a consequence of Huisken's type monotonicity formula and
\eqref{e4i}, there holds
\begin{equation}\label{e7i}
\frac{d}{dt}\int_{\re^N} \psi(x,t) d\mu^\e_t(x)\le 0\,,
\end{equation}
$i.e.$, monotonicity of the function $t\mapsto \int_{\re^N}
\psi(x,t) d\mu^\e_t(x)$; here, for each fixed $y\in \re^N, s>0$,
\begin{equation}\label{e9i}
\psi(x,t)\equiv \psi(x,t; y,s)
:=\frac{e^{-\frac{|x-y|^2}{4(s-t)}}}{[4\pi(s-t)]^{\frac{N-1}2}}
\quad \textrm{for all}\;\; x\in \re^N, 0\leq t<s\,;
\end{equation}
observe that this function $\psi$ is, up to a multiplicative factor $\sqrt{4\pi(s-t)}$, exactly the {\it
backward heat kernel} in dimension $N$.

Next it is shown that then there are a Radon measure $\mu_t$ on
$\mathbb{R}^N$ and a sequence $\{\e_n\}\subset (0, \infty),\, \e_n\to 0$ as
$n\to\infty$ such that, for every $t>0$, $\mu^{\e_n}_t$ converges
as Radon measure on $\mathbb{R}^N$ to $\mu_t$ for all $t\ge 0$ as $n\to
\infty$. These measures are shown to be
$(N-1)-$rectifiable, as a consequence of density bounds derived
from \eqref{e7i}. Finally, Brakke's inequality for $d\mu_t$ is obtained
from the corresponding approximate ones valid for $d\mu_t^\e$.

\medskip

The aim of this paper and of \cite{PiPu2} is to generalize the
results in \cite{Ilm1} recalled above, to the case of solutions
$u^\e$ to problem \eqref{e1}-\eqref{e1a} on Riemannian manifolds. We always assume that
there exists $\l\in \re$ such that
\begin{equation}\label{e53}
\Ric(X,X)\geq \l \langle X, X\rangle\quad\textrm{for all}\;\; X\in
T_x M, x\in M\,;
\end{equation}
here $\Ric$ denotes the Ricci tensor on $M$. Note that for
$M=\re^N$, we have $\l=0$; for the hyperbolic space $\mathbb H^N$,
$\l=-(N-1)$; for the sphere $\mathbb S^N$, $\l=N-1$  (see
Subsection \ref{DG}). Indeed, in these cases \eqref{e53} holds with the equality sign.

Observe that, under the assumption on the Ricci curvature, we can
apply comparison principle for \eqref{e1}-\eqref{e1a}. In addition, we can
treat not only well prepared initial conditions, but also quite general
initial conditions. Hence our results with $M=\re^N$ extend those
in \cite{Ilm1} in this respect.

Note that while mean curvature flow has been investigated also on
Riemannian manifolds (see, $e.g.$ \cite{AAM}, \cite{Hui1},
\cite{Hui2}, \cite{Ilm3}), to the best of our knowledge, the
question of approximation of mean curvature flow via Allen-Cahn
equation on Riemannian manifolds has not been addressed. On the other hand, the connection between the
stationary Allen-Cahn equation and minimal hypersurfaces has been widely studied e.g.
in
\cite{PacRit}, \cite{dPKWY} and \cite{PisPons}.

Now, we outline results that will be shown in the present paper
and we briefly mention the content of \cite{PiPu2}. For any
$\e>0$, define the energy density
\begin{equation}\label{e8i}
d\mu^\e_t := \left\{\frac{\e}2 |\nabla u^\e|^2 + \frac
1{\e}F(u^\e)\right\}d\mathcal V(x),
\end{equation}
$u^\e$ being a solution to equation \eqref{e1} and $d\mathcal V$
the volume element on $M$. For $d\mu^\e_t$ we shall prove {\it
Brakke's} and {\it Huisken's} type formulas (see Lemma
\ref{lemma1}, and respectively, Lemma \ref{lemma5}). Also in this
case, they contain the discrepancy Radon measure
\[d \xi_t^\e:=\left(\frac{\e}2 |\nabla u^\e|^2-\frac 1{\e} F(u^\e)\right)d \mathcal V(x)\,.\]

For a general class of initial conditions, without supposing that
$u_0^\e$ are somehow well prepared, we prove that
\begin{equation}\label{u13}
\limsup_{\e\to 0^+} \sup_{(x,t)\in Q} \xi^\e_t(x)\leq 0\,,
\end{equation}
for each compact subset $Q\subset M\times (0,\infty)$. To do this,
we adapt and improve an elementary but very clever idea from \cite{HT}. Observe that in
\cite{HT} the stationary problem for $M=\re^N$ is addressed.
Moreover, it is only shown that the positive part of $\xi^\e$ is
bounded, uniformly with respect to $\e$. In order to show
\eqref{u13} we improve some estimates in the argument of \cite{HT}
and extend them to the case of Riemannian manifolds (see Section
\ref{gic}).

\smallskip

However, for well prepared initial conditions $u_0^\e$ also an
alternative strategy can be used. In fact, in Section \ref{wpic}
we prove for properly well prepared initial conditions an
asymptotic control of discrepancy, by methods similar to those
used in \cite{Ilm1}. However, some differences from \cite{Ilm1}
occur, for the presence of the general Riemannian metric on $M$,
which we describe below.

As a consequence, instead of $q_\e$, it is convenient to consider
the one-dimensional profile $h^\e$, which is the solution, for any
$\e>0$, to problem
\begin{equation}\label{e38}
\left\{
\begin{array}{ll}
\frac 1{\varphi}\big\{\varphi h_{\e}' \big\}'= f_\e(h_{\e})
 & \textrm{in} \;\;  (0,1)
\\& \\
h_\e(0)=0,\;\; h_\e(1)=1\,,
\end{array}
\right.
\end{equation}
(see Subsection \ref{1d}) where $\varphi:[0,1]\to(0,\infty)$ is an
increasing convex smooth function such that $\varphi(0)\ge 1,
\varphi'(0)=0$ that will be chosen to balance some curvature effects.
 We still denote by $h_\e$ the odd reflection of
the solution of \eqref{e38}. Indeed, note that the ordinary
differential equation in \eqref{e38}, for the choice
$\varphi\equiv 1$, coincides with that solved by $q_\e$ but for technical reasons
it is more convenient to consider \eqref{e38} on a bounded interval. As a
preliminary step we shall prove that (see Subsection \ref{1d})
\begin{equation}\label{e161}
\limsup_{\e \to 0} \sup_{\t\in (0,1)}\e \left\{\frac 1
2h_\e'(\t)\,^2- F_\e(h_\e(\t))\right\}\, \le\, 0\,,
\end{equation}
exploiting the fact that \eqref{e38} is now solved in a bounded
interval.

Concerning the initial conditions $u^\e_0$, well prepared data
will be now of the form  (see Subsection \ref{ic}):
\begin{equation}
\label{wellprepdata}
u_0^{\e}(x):= h^\e\left( \Psi ( \tilde d(x, \Sigma_0) )
\right)\, \quad (x\in M)\,,
\end{equation}
where $E_0\subset M$ is an open bounded subset with smooth
boundary $\Sigma_0:=\pa E_0$. Now, $d(x, \Sigma_0)$  is the
Riemannian distance of $x\in M$ to $\Sigma_0$, while the signed
distance $\tilde d(x,\Sigma_0)$ is defined accordingly and $\Psi(\cdot)$ is a suitable smoothed and $1-$Lipschitz truncation of the identity, which makes $u_0^\e$ constant far from $\Sigma_0$, where the distance function is possibly singular.

\smallskip

Define the function $z^\e:M\times[0,\infty)\to \re$ by
\begin{equation}\label{e71}
u^\e(x,t):= h_\e \big(z^\e(x,t) \big)\quad \big(x\in M, t\ge
0\big)\,.
\end{equation}
Under the assumption
\begin{equation}\label{e55b}
|\nabla z^\e(x,0)|\le 1\quad \textrm{for all}\;\;\, x\in M\,,
\end{equation}
 which clearly follows if \eqref{wellprepdata} holds, and
\begin{equation}\label{e51}
\left(\frac{\varphi'}{\varphi}\right)'\ge \max\{-\l, 0\}\quad
\textrm{in}\;\;(0,1)\,,
\end{equation}
we shall prove that (see Subsection \ref{discrep})
\begin{equation}\label{e52}
|\nabla z^\e(x,t)|\le 1\quad \textrm{for all}\;\; x\in M, t\ge
0\,.
\end{equation}
Indeed, note that, in view of \eqref{e55b} and \eqref{e51}, we can infer
that $\tilde w\equiv 1$ is a supersolution to a certain parabolic
equation solved by $|\nabla z^\e|^2$. Hence, from maximum
principle we can prove that \eqref{e52} holds true. This argument
fails, if  $\varphi\equiv 1$ and $\lambda<0$, because of the presence of
an extra term related to the Ricci curvature; so, we cannot
consider $q^\e$ instead of $h^\e$, e.g. when $M=\mathbb H^N$, and this is the reason why we have to introduce the profile $h_\e$ defined in \eqref{e38}. As a
consequence of \eqref{e161} and \eqref{e52} we have \eqref{e161}.

%\begin{equation}\label{e55}
%\limsup_{\e\to 0^+} \e\left\{\frac 1 2|\nabla u^\e|^2- \frac
%1{\e^2} F(u^\e)\right\}\,d\mathcal V(x) \le 0\,.
%\end{equation}

\medskip

Note that both for general initial conditions and for well
prepared initial conditions we cannot prove that the discrepancy
term is nonpositive, as occurred in \eqref{e4i} in the Euclidean space. However,
condition \eqref{u13} will play the same role as \eqref{e4i} has
in the case $M=\re^N$.

Then, in Section 5 from Huisken's type equality for the density energy and
\eqref{u13} we obtain the following inequality (see Theorem
\ref{prop3})
\begin{equation}\label{e60}
\begin{split}
\frac{d}{dt}\int_{M} \phi(x,t)d\mu^\e_t\,\le\,
\frac{C_3}{\sqrt{s-t}}\int_{M} \phi(x,t)d\mu^\e_t + C_4 + \frac
{C_5}{\sqrt{s-t}}
\end{split}
\end{equation}
for all $0\le t<s$, for some positive constants $C_3, C_4, C_5$
independent of $\e$. Inequality \eqref{e60} is a natural counterpart on a manifold of
the monotonicity formula \eqref{e7i} but,  due to the presence of extra terms,
\eqref{e60} does not imply monotonicity for the function
$t\to\int_{M} \phi(x,t)d\mu^\e_t(x)$.

Here, for any fixed reference point $(y,s)\in
M\times (0,\infty)$, $\phi(x,t)\equiv \phi(x,t; y,s)$ is a
suitable kernel, which replaces \eqref{e9i}. It depends explicitely on the
Riemannian distance $d(x)=d(x,y)$ for $x,y\in M$ as follows
\begin{equation}
\label{kernel}
\phi(x,t)=\hat{\zeta}(d^2(x)) (s-t)^{-\frac{N-1}{2}} e^{-\frac{d^2(x)}{4(s-t)}} \, ,
\end{equation}
 furthermore,
in constrast with the case of $\re^N$, it has a suitably small
compact support in space due to the cut-off function $\hat{\zeta}$.
 As a consequence,
it allows us to
control the behavior of $d\mu^\e_t$ only at small scales. For this
reasons, we shall refer to \eqref{e60} as a {\it local almost
monotonicity formula}. This choice of the kernel is very natural, since, up to the cut-off and the factor $\sqrt{s-t}$, is nothing but the leading order term in the expansion of the backward heat kernel on the manifold $M$ with pole at $(y,s)$ for short times.
It would be very interesting to find a more precise localized monotonicity formula for the Allen-Cahn equation on a manifold containing no error term. It should be analouge to the one in \cite{Ilm3} for $\mathbb{R}^N$ but local as the celebrated formula in \cite{Ek} for the mean curvarure flow, still in $\mathbb{R}^N$.

As a consequence of \eqref{e60} we obtain, for all $0\leq
t_0<t<s$,
\begin{equation}\label{e61}
\begin{split}
\mathcal G(t)\le e^{\frac{C_3}2 (\sqrt{s-t_0}  -\sqrt{s-t})}\big[ \mathcal
G(t_0) + C_4(t-t_0) + C_5( \sqrt{s-t_0}  -\sqrt{s-t}) \big]\,,
\end{split}
\end{equation}
where
\[\mathcal G(t):= \int_{M} \phi(x,t) d\mu^\e_t\quad (0\le t<s)\, \]

\smallskip

and this is precisely the inequality needed to have uniform density bounds for the measures $\mu^\e_t$ at small scales. We conclude Section 5 giving some useful compactness properties for the solutions $u^\e$ both in $L^1_{loc}$ and in the space of functions of bounded variation.

Finally, let us mention that, out of its independent interest,
inequality \eqref{e61} will be used in \cite{PiPu2} to prove that
there exist a Radon measure $\mu_t$ on $M$ and a sequence
$\{\e_n\}\subset (0, \infty),\, \e_n\to 0$ as $n\to\infty$ such
that, for every $t>0$, $\mu^{\e_n}_t$ converges as Radon measure
on $M$ to $\mu_t$ for all $t\ge 0$ as $n\to \infty$. Moreover,
$\mu_t$ will be $(N-1)-$rectifiable and they will satisfy the
Brakke's inequality, i.e. they will be a generalized solution of the mean curvature flow in the sense of varifolds with the surface measure on $\Sigma_0$ as initial data.

\section{Preliminaries from Differential Geometry}\setcounter{equation}{0}\label{DG} In this Section we
recall some basic facts and notations from Riemannian Geometry,
that will be used in the sequel and in \cite{PiPu2}, too (for more
details see, $e.g.$, \cite{GHL}, \cite{Jost}).

\smallskip
Let $M$ be an $N-$dimensional Riemannian manifold, equipped with a
metric tensor $g$. For any given point $x\in M$, let $T_x M$ be
the tangent space at $x$, $TM$ be the tangent bundle, $T_x^*M$ be
the cotangent space at $x$, $T^*M$ be the cotangent bundle,
$\Gamma(TM)$ denote the vector space of smooth sections of $TM$,
$i.e.$ the smooth vector fields on $M$. In local coordinates
$\{x^1,\ldots, x^N\}$, we have a natural local basis
$\left\{\frac{\pa}{\pa x^1},\ldots, \frac{\pa}{\pa x^N}\right\}$
for $TM$. The metric tensor $g=g_{ij} dx^i\otimes dx^j$ is
represented by a smooth matix-valued function
$g_{ij}=g\left(\frac{\pa}{\pa x^i}, \frac{\pa}{\pa x^j}\right)$,
so that locally the {\it inner Riemannian product} $\langle\cdot,
\cdot\rangle$ is given by
\begin{equation}\label{e116}
\langle X, Y \rangle :=  g_{ij}\, X^i Y^j\,,
\end{equation}
where the vectors $X=X^i\frac{\pa}{\pa x^i}$, $Y=Y^i\frac{\pa}{\pa
x^i}$ belong to the tangent space $T_{x}(M)$. The induced {\it
geodesic distance} between any two points any $x,y\in M$ will be
indicated by $d(x,y)$. For any $x_0\in M, r>0$ let
$B_r(x_0):=\big\{x\in M\,|\,d(x,x_0)<r\big\}\,.$ The {\it
gradient} $\nabla u$ of a function $u\in C^1(M)$ is given by
\[
\big(\nabla u\big)^i := g^{i j}\frac{\pa u}{\pa x^j} \quad \qquad
(i=1,\ldots, N)\,,
\]
so that
\[ du(X)\,=\, \langle X ,\nabla u\rangle\quad \textrm{for any}\;\;
X\in \Gamma(TM)\,.
\]

Recall that the for any vector field $Y\in \Gamma(TM)$ there
exists a unique smooth function on $M$, denoted by $\di Y$, such
that the following identity holds:
\[\int_{M} \phi\, \di Y\, d\mathcal V\,=\,-\int_{M} \langle Y, \nabla \phi\rangle\, d \mathcal V\,\]
for all $\phi\in C^1_c(M)$. Furthermore, in local coordinates
\[\di Y\,=\, \frac{1}{\sqrt{\bar g}}\frac{\pa}{\pa x^k}\left(\sqrt{\bar g}Y^k\right)\,,\]
where $\bar g:=\,\hbox{det}\, (g_{ij})$.

The {\it Laplace-Beltrami} operator on $M$ is given by:
\[\Delta = \di\circ \nabla = \frac{1}{\sqrt{\bar g}}\frac{\pa}{\pa x^i}\left(\sqrt{\bar g}g^{ij}\frac{\pa}{\pa x^j} \right)\,.\]

\smallskip

The {\it Levi-Civita} connection $\cd$ of the metric $g$ is given
by
\[\cd_{\frac{\pa}{\pa x^i}}{\frac{\pa}{\pa x^j}}=\Gamma^k_{ij}\frac{\pa}{\pa x^k},\]
where
\begin{equation}\label{e118}
\Gamma^k_{ij}:=\frac 1 2g^{kl}\left(\frac{\pa g_{jl}}{\pa x^i}
+\frac{\pa g_{il}}{\pa x^j}-\frac{\pa g_{ij}}{\pa x^l}\right)
\end{equation}
are the {\it Christoffel} symbols.

\smallskip

We also recall that the Hessian of $f\in C^2(M; \re)$ is the
symmetric endomorphism of $TM$ defined by
\[\Hess f(X):= \cd_X \nabla f\,\quad \textrm{for any}\;\; X\in \Gamma(TM)\,,\]
or its associated symmetric bilinear form on $TM$ defined by
\[(\Hess f)(X,Y):=X(Y(f))- \cd_XY(f)\quad\textrm{for any}\;\; X,Y\in \Gamma(TM)\,.\]
We have:
\[(\Hess f)(X,Y)=\langle \cd_X(\nabla f) , Y \rangle \quad \textrm{for any}\;\; X,Y\in \Gamma(TM)\,.\]
Also, in local coordinates, there holds:
\[(\Hess f)(X,Y)=\sum_{i,j=1,N} X^i X^j (\Hess f)_{ij},\]
where
\begin{equation}\label{e119}
(\Hess f)_{ij}=(\Hess f)\left(\frac{\pa}{\pa x^i}, \frac{\pa}{\pa
x^j}\right)=\frac{\pa ^2f}{\pa x^i\pa x^j} -\Gamma^k_{ij}\frac{\pa
f}{\pa x^k}\,.
\end{equation}

\smallskip

In terms of the Hessian, the Laplace-Beltrami operator rewrites
as:
$$\Delta f = \sum_{i,j=1}^N g^{ij}(\Hess
f)\left(\frac{\pa}{\pa x^i}, \frac{\pa}{\pa x^j}\right)= tr\,(\Hess f)\,;$$ here
and hereafter $\textrm{tr}$ denotes the trace operator (taken
fiberwise).

\smallskip
For any $y\in M$, denote by $\inj(y)$ the injectivity radius at
$y$. In the sequel we use the next lemma.
\begin{lemma}\label{lemma1u}
Let $y\in M, d(x):=d(x,y)$ for all $x\in M$. For any compact
subset $K\subset M$ there exists a constant $C>0$ such that if
\[y\in K,\;\; d(x)\leq \frac 12 \inf_{y\in K} \inj (y),  \]
then
\begin{equation}\label{u50}
\left \|\frac 12 \Hess d^2(x)(X,X)- g(X,X)\right \| \leq C
d^2(x)\| X\|^2 \quad \textrm{for any}\;\; X\in T_xM\,,
\end{equation}
and
\begin{equation}\label{u51}
\left| \frac 1 2 \Delta d^2(x) - N    \right|\leq C d^2(x)\,.
\end{equation}
\end{lemma}

\medskip

For $u, v, \phi\in C^2(M; \re)$, it is direct to see that:
\begin{equation}\label{e5}
\langle\nabla u, \nabla\langle\nabla \phi,\nabla v \rangle \rangle
= (\Hess \phi)(\nabla u, \nabla v)\,+\, (\Hess v)(\nabla u, \nabla
\phi)\,.
\end{equation}

The {\it curvature tensor} of the Levi-Civita connection $\cd$ is
given by
\[R(X,Y)Z:= \cd_X \cd_Y Z- \cd_Y \cd_X Z - \cd_{[X,Y]}Z\quad \textrm{for any}\;\, X,Y,Z\in \Gamma(TM)\,;\]
in local coordinates,
\[R\left(\frac{\pa}{\pa x^i}, \frac{\pa}{\pa x^j}\right)\frac{\pa}{\pa x^l}=R^k_{lij}\frac{\pa}{\pa x^k},\]
where
\[R^k_{lij}:= \frac{\pa \Gamma^k_{jl}}{\pa x^i}- \frac{\pa \Gamma^k_{il}}{\pa x^j}+\Gamma^k_{im}\Gamma^m_{jl}-\Gamma^k_{jm}\Gamma^m_{il}\,.\]

The sectional curvature of the plane $X \wedge Y$ spanned by the
linearly independent tangent vectors $X=X^i\frac{\pa}{\pa x^i},
Y=Y^i\frac{\pa}{\pa x^i}\in T_xM$ is
\[K(X\wedge Y):=\frac{\langle R(X,Y)X  ,Y  \rangle}{|X\wedge Y|^2},\]
where $|X\wedge Y|^2= |X|^2 |Y|^2-\langle X,Y\rangle^2\,.$ The
{\it Ricci curvature} in the direction $X=X^i\frac{\pa}{\pa
x^i}\in T_xM$ is
\[ \Ric(X,X):= g^{jl}\langle R\left(X, \frac{\pa}{\pa x^j}\right)X, \frac{\pa}{\pa x^l}\rangle;\]
the {\it Ricci tensor} is
\[ R_{ik}=g^{jl}R_{ijkl}=R_{ki},\]
where $R_{ijkl}=g_{im}R^m_{jkl}\,.$

Furthermore, recall the {\it Bochner-Weitzenb$\ddot{o}$ch formula}
: for $\phi\in C^2(M;\re)$ there holds
\begin{equation}\label{e53a}
\frac 1 2 \Delta(|\nabla \phi|^2)\,=\, |\Hess \phi|^2 + \langle
\nabla \phi, \nabla\Delta\phi \rangle + \Ric(\nabla \phi,
\nabla\phi)\,.
\end{equation}

\smallskip

\section{Asymptotic control of discrepancy for general initial
conditions}\setcounter{equation}{0}\label{gic}
\subsection{General initial conditions}
Let $E_0\subset M$ be an open bounded subset with $C^2-$boundary
$\pa E_0=\Sigma_0$. Note that there exist $R_0>0, C_0>0$ such that
$$\mathcal H^{N-1}\big(\Sigma_0\cap B_R(x)\big)\le C_0
\omega_{N-1}R^{N-1}$$ for all $0<R<R_0$\,.

\smallskip

For any $\e>0$ set
\begin{equation}\label{e2}
E^\e(x,t):=\frac 1 2 |\nabla u^\e|^2+\frac 1 {\e^2} F(u^\e)\quad
(x\in M, t\ge0)\,;
\end{equation}
clearly (see \eqref{e8i}),
\begin{equation}\label{e3}
d\mu^\e_t(x)=\e E^\e(x,t) d \mathcal V(x)\quad (x\in M, t\ge 0)\,.
\end{equation}

\smallskip
Concerning the initial conditions $u_0^\e$ (and the corresponding
$\mu^\e_0\equiv \mu^\e(\cdot,0)$ given by \eqref{e8i}) we always
assume the following:
\[\left\{
\begin{array}{l}
(i)\,\,\mu^\e_0\to \a  \mathcal H^{N-1}\lfloor
\Sigma_0\,\;\hbox{as}\,\, \e\to 0\,\,\hbox{as Radon measures, for
some}\,\,\a\ge 0\,;
\vspace{.2 cm} \\
(ii)\,u_0^\e \to 2 \chi_{E_0} -1\quad \textrm{as}\;\,\e\to\infty\,\,\textrm{in}\;\, BV_{loc}(M)\, \textrm{weakly}-*; \\
(iii)\,\hbox{there exists}\,\, C_0>0\,\,\hbox{such that}\,\,
\frac{\mu_0^\e(B_R(x))}{\omega_{N-1}R^{N-1}}\le C_0 \\ \quad
\;\;\;\;
\textrm{for all}\; x\in M,\, 0<R<R_0,\, 0<\e<1; \vspace{.2 cm} \\
(iv)\, \hbox{there exists}\,\, k_0>0\,\, \hbox{such that}\,\,
\|u_0^\epsilon\|_\infty\leq k_0\,;\\
(v)\,\; u_0^\e \in C^1(M) \hbox{ and there exists}\,\, \check
C>0\,\,\hbox{such that for any}\\ \quad \;\,\;\, 0<\e<1\; \|\nabla
u_0^\e\|_\infty\le
\frac{\check C}{\e} \, . \vspace{.2 cm} \\
\end{array}
\right. \leqno(H_1)
\]

Throughout this section, we will not assume any further structure assumption on the initial data and, on the contrary, even the previous hypoteses both on $\Sigma_0$ and on $u^\e_0$ could be further relaxed.

\subsection{Global existence and uniqueness results}
Concerning existence and uniqueness of solutions to problem
\eqref{e1}-\eqref{e1a} we state the next Proposition.

\begin{proposition}\label{prop11}
Let hypotheses $(H_0), (H_1)$ be satisfied. Then problem
\eqref{e1}-\eqref{e1a} admits a unique bounded solution. Moreover,
$u^\e\in C^{\infty}\big(M\times(0,\infty)\big)\cap
C^0\big(M\times[0,\infty)\big),$ and
\begin{equation}\label{e121b}
|u^\e | \leq  k_0 \quad \textrm{for all}\;\; x\in M, t>0\,.
\end{equation}
In addition,
\begin{equation}\label{e162}
\sup_{\e>0} \sup_{t\in (0,\infty)} \e \int_M E^\e d\mathcal V \le
C_2\,
\end{equation}
where $C_2:=\sup_{\e>0}\mu_0^\e(M)$;
\begin{equation}\label{e150}
t\mapsto \int_M E^\e(x,t) d\mathcal V(x)\,\;\,\textrm{is
nonincreasing for}\;\; t>0\,.
\end{equation}
\end{proposition}
\noindent{\it Proof\,.} Existence and regularity of solutions can
be shown by usual methods, e.g. solving the corresponding IBV problems on an increasing family of bounded domains with smooth boundary and arguing by local a-priori estimates and compactness. In view of \eqref{e53}, from results in
\cite{Dod} uniqueness and comparison principles for problem
\eqref{e1}-\eqref{e1a} can be easily deduced. In view of
$(H_0)-(ii)$, the functions $\bar v\equiv k_0, \underline v\equiv
- k_0$ are a supersolution and, respectively, a subsolution to
problem \eqref{e1}, \eqref{e1a}. Hence, by comparison principle
\eqref{e121b} follows. Finally, inequality \eqref{e162} and the
property \eqref{e150} follows passing to the limit in the global energy inequality on the approximating domains.
\hfill $\square$

\begin{proposition}\label{prop11b}
Let hypotheses $(H_0), (H_1)-(iv)$ be satisfied. Let $u^\e$ be the solution to problem \eqref{e1}-\eqref{e1a}. Then \eqref{e121b} holds true.
Furthermore, for any compact subset $K\subset M$ and for any $\t\in (0,T)$ there exists a constant $\tilde k>0$ such that
\begin{equation}\label{e163}
\|\nabla u^\e(\cdot, t)\|_{L^\infty(K)}\le \frac{\tilde
k}{\e}\quad \textrm{for all}\;\; t\in (\t, T)\,;
\end{equation}
\begin{equation}\label{e163b}
\e \xi^\e(x,t) \leq \tilde k \quad \textrm{for all}\;\; x\in K,
t\in (\t, T).
\end{equation}
\end{proposition}

\noindent{\it Proof\,.} Note that \eqref{e121b} can be deduced as in the proof of Proposition \ref{prop11}. Moreover, \eqref{e163} follows by standard parabolic estimates, writing the equation in local coordinates and arguing by scaling. Consequently,\eqref{e163b} is obtained, in view of \eqref{e121b}.
The proof is complete. \hfill $\square$

\begin{remark}\label{remp11} For further references, note that from $(H_1)-(i)$ and \eqref{e150} it is
direct to see that, for each compact subset $K\subset M, T>0,
\t\in [0, T)$, there holds:
\begin{equation}\label{e170}
\sup_{\e>0} \sup_{t\in (\t,T)} \e \int_K E^\e d\mathcal V \le
\overline C\,
\end{equation}
for some constant $\overline C>0$ depending on the compact subset
$K, \t>0, T>0$, and independent of $\e>0$. Indeed, under $(H_1)-(i)$ we have
\begin{equation}\label{e170b}
\overline C\leq C_2 \, ,
\end{equation}
where $C_2$ given in Proposition \ref{prop11} is clearly independent of $K, \t>0, T>0$ in view of \eqref{e162} and \eqref{e150}. However, in
the sequel most of the time the arguments will rely only on \eqref{e170} and we shall not use the
property \eqref{e170b}.
\end{remark}

\subsection{Asymptotic control of discrepancy}
We prove the next result.
\begin{proposition}\label{propu1}
Let assumption $(H_0)$ hold true. Let $\{u^\e\}$ be a family of uniformly bounded
solutions to problem \eqref{e1}-\eqref{e1a}, i.e. \eqref{e121b} is verified. Then \eqref{u13} is
satisfied.
\end{proposition}

In order to prove Proposition \ref{propu1} we need some
preliminary results.

\begin{lemma}\label{lemmau1}
Let $\bar x\in M, r>0, \bar t>4 r^2.$ Let $\Omega_0=B_{4r}(\bar
x)\times (\bar t- 4r^2, \bar t], \Omega:= B_r(\bar x)\times (\bar
t-r^2, \bar t].$ Assume that, for some $C>0$,
\begin{equation}\label{u1}
\sup_{0<\e<1}\|u^\e\|_{L^\infty(\Omega_0)}\leq C.
\end{equation}
Then, for any $\s_0\in (0,2)$, there exists a constant
$C_0=C_0(\Omega, C, \s_0)>0$ such that
\begin{equation}\label{u2}
\|u^\e\|_{L^\infty(\Omega)}\leq 1+ C_0 \e^{\s_0}\quad \textrm{for
any}\;\; \e \in (0,1)\,.
\end{equation}
\end{lemma}

\noindent{\it Proof\,.\,\,} Fix any $\s_0\in (0,2)$. It suffices
to show the thesis with $C_0=1$ and $\e\to 0^+$. In fact, as a
consequence of this, we can immediately get \eqref{u2}, taking
possibly a bigger $C_0$.

\smallskip

Suppose, by contradiction, that there exists a sequence
$\{\e_n\}\subset (0,1)$ such that $\e_n\to 0^+$ as $n\to \infty$
and
\[\sup_{\Omega} u^{\e_n}\geq 1+ \e_n^{\s_0}\,.\]
The case $\inf_{\Omega} u^{\e_n}\geq -1-\e_n^{\s_0}$ can be
treated with obvious modifications; so we do not discuss it in
details.

Let $\Omega_1:= B_{2r}(\bar x)\times (\bar t-2r^2, \bar t]$. For
any $n\in \ene$ select $\varphi_n\in C^\infty(\bar \Omega_0)$ such
that
\[\varphi_n\equiv 1 + \frac 1 2\e_n^{\s_0}\quad \textrm{in}\;\; \Omega,\]
\[1< 1+\frac 1 2 \e_n^{\s_0}\leq \varphi_n\leq 1 + C\quad \textrm{in}\;\; \Omega_0,\]
\[\varphi_n\equiv 1+ C\quad \textrm{in}\;\; \Omega_0\setminus \Omega_1\,;\]
moreover, for some $\bar C>0$, for all $n\in \ene$,
\begin{equation}\label{u3}
|\nabla \varphi_n|+ |\Hess (\varphi_n)|+ |\pa_t \varphi_n|\leq
\bar C\quad \textrm{in}\;\; \Omega_0\,.
\end{equation}

Set
\[g_n:= u^{\e_n}-\varphi_n,\]
so that
\[g_n\leq -1\quad \textrm{in}\;\; \big[B_{4r}(\bar x)\times\{\bar t-4r^2\}\big]\cup\big[\pa B_{4r}(\bar x)\times (\bar t-4r^2, \bar t]\big]\,.\]
Furthermore,
\[\sup_{\Omega_0} g_n \geq \sup_{\Omega} g_n \geq \frac 1 2\e_n^{\s_0}>0\,.\]
Then $\max_{\bar \Omega_0}g_n=g_n(x_n, t_n)$ for some $(x_n,
t_n)\in \Omega_0\,.$ Thus, using \eqref{e1}, the fact that
$u^{\e_n}(x_n, t_n)>1, (H_0)-(iv)$ and \eqref{u3} we obtain
\[0\geq \Delta g(x_n, t_n)-\pa_t g(x_n, t_n)= \Delta u^{\e_n}(x_n, t_n)-\pa_t u^{\e_n} (x_n, t_n)\]\[+\pa_t\varphi_n(x_n, t_n)-
\Delta \varphi_n(x_n, t_n)\]
\[=\frac{F'(u^{\e_n}(x_n, t_n))}{\e_n^2}+\pa_t \varphi_n(x_n, t_n)-\Delta \varphi_n(x_n, t_n)  \]
\[\geq \frac{F'(u^{\e_n}(x_n, t_n))}{\e_n^2}-\frac{F'(\varphi_n(x_n, t_n))}{\e_n^2}+ \pa_t \varphi_n(x_n, t_n)-\Delta \varphi_n(x_n, t_n)\geq
\]
\[=\frac{F'((1-s)\varphi_n+su^{\e_n})}{\e_n^2}(x_n, t_n)\Big|_{s=0}^{s=1}+\pa_t \varphi_n(x_n, t_n)-\Delta \varphi_n(x_n, t_n)
\]
\[=\int_0^1\frac{d}{ds}\frac{F'((1-s)\varphi_n+su^{\e_n})}{\e_n^2}(x_n, t_n)ds +\pa_t \varphi_n(x_n, t_n)-\Delta \varphi_n(x_n, t_n) \]
\[=\frac{g_n(x_n, t_n)}{\e_n^2}\int_0^1 F''((1-s)\varphi_n+su^{\e_n})(x_n, t_n)ds+\pa_t \varphi_n(x_n, t_n)-
\Delta \varphi_n(x_n, t_n) \]
\[\geq \frac 12 \inf_{1<s<1+C}F''(s)\e_n^{\s_0-2}-\bar C.\]
This is clearly impossible for $n\in \ene$ large enough, hence the
thesis follows. \hfill $\square$

\medskip

Define
\[G(u):= \e^{\sigma}(2 H_0 - u^2),\]
so that
\[G>0, \quad G''=-2\e^\s <0\, \]
Set
\[\phi_{G}^{\e}:= \e \xi_t^\e - G=\frac{\e^2}2|\nabla u^\e|^2 - F(u^\e)- G(u^\e)\,,\]
where $u^\e$ is a solution to equation \eqref{e1}. Hence, we have:
\begin{equation}\label{u4}
\begin{split}
(\Delta -\pa_t)\phi_{G}^{\e}= (\Delta -\pa_t)\frac{\e^2}2|\nabla
u^\e|^2- (\Delta -\pa_t)(F+G)\hspace{1.2 cm}\\
=\e^2 \big[|\Hess u^\e(\nabla u^\e, \nabla u^\e)|^2   + \langle
\nabla u^\e, \nabla \Delta u^\e\rangle + \Ric(\nabla u^\e, \nabla
u^\e)\hspace{1 cm}\\
-\langle\nabla u^\e, \nabla \pa_t u^\e \rangle \big]+ (F'+G')\pa_t
u^\e - \di\big((F'+G')\nabla u^\e\big)\hspace{.8 cm}
\\=\e^2\big[|\Hess u^\e(\nabla u^\e, \nabla u^\e)|^2+ \langle \nabla u^\e ,\nabla (\Delta
u^\e-\pa_t u^\e) \rangle +\Ric(\nabla u^\e, \nabla u^\e)\big]\\
+(F'+G')(\pa_t u^\e-\Delta u^\e)-(F''+G'')|\nabla u^\e|^2 \hspace{2.2 cm}\\
=\e^2|\Hess u^\e(\nabla u^\e, \nabla u^\e)|^2 + \e^2\Ric(\nabla
u^\e, \nabla u^\e)\hspace{2.2 cm}\\- \frac{1}{\e^2}(F'+G')F' -
G''|\nabla u^\e|^2\,.\hspace{2.8 cm}
\end{split}
\end{equation}
Note that
\begin{equation}\label{u5}
|\nabla u^\e|^2|\Hess u^\e(\nabla u^\e, \nabla u^\e)|^2 \geq \frac
1 2 \big|\nabla|\nabla u^\e|^2\big|^2\,.
\end{equation}
To see this, take any $p\in M$ and fix an orthonormal frame
$\{E_i\}_{i=1,\ldots, N}$ around $p$. Thus,
\[|\nabla u^\e|^2|\Hess u^\e (\nabla u^\e, \nabla u^\e)|^2 =\sum_{i=1}^N|\nabla u^\e|^2|\cd_{E_i} \nabla u^\e|^2  \]
\[\geq \sum_{i=1}^N\big(\langle \nabla u^\e ,\cd_{E_i}\nabla u^\e \rangle \big)^2=\frac 1 2\sum_{i=1}^N\big(E_i|\nabla u^\e|^2\big)^2
=\frac 1 2 \Big|\nabla |\nabla u^\e|^2 \Big|^2\,.\] So, \eqref{u5}
has been verified. From \eqref{u4}, \eqref{u5} and \eqref{e53} we
deduce that, whenever $\nabla u^\e\neq 0$,
\[\e^2 |\Hess u^\e(\nabla u^\e, \nabla u^\e)|^2 \geq \frac 1{\e^2 |\nabla u^\e|^2}\big|\nabla (\phi^\e_G+F+G) \big|^2 \]
\[\geq 2\frac{F'+G'}{\e^2|\nabla u^\e|}\frac{\nabla u^\e}{|\nabla u^\e|}\nabla \phi_G^\e + \frac 1{\e^2}(F'+G')^2,\]
therefore,
\[(\Delta -\pa_t)\phi_G^\e - 2 \frac{F'+G'}{\e^2|\nabla u^\e|}\left\langle \frac{\nabla u^\e}{|\nabla u^\e|},\nabla \phi_G^\e\right\rangle\]
\[\geq \frac 1{\e^2}(G')^2+ \frac 1{\e^2}F'G' + \l |\nabla u^\e|^2\e^2 - G''|\nabla u^\e|^2\,.\]

We summarize these computations in the following result.
\begin{lemma}\label{lemmau2}
Whenever $\nabla u^\e\neq 0$, let
\begin{equation}\label{u15}
\mathcal A^\e:=(\Delta - \pa_t)\phi_G^\e +
\left(\frac{2G''}{\e^2}-2\l \right)\phi_G^\e
-\frac2{\e^2}\frac{(F'+G')} {|\nabla u^\e|}\left\langle
\frac{\nabla u^\e}{|\nabla u^\e|},\nabla \phi_G^\e\right\rangle;
\end{equation}
\begin{equation}\label{u16}
\mathcal B^\e:=\left(2\l - \frac 2{\e^2}G''\right)(F+G)+\frac
1{\e^2}(G')^2+ \frac{F'G'}{\e^2}\,.
\end{equation}
Then there holds:
\begin{equation}\label{e17}
\mathcal A^\e \geq \mathcal B^\e\,.
\end{equation}
\end{lemma}
The next result is an improvement of \cite{HT}.

\begin{lemma}\label{lemmau3}
Let $\bar x\in M, r>0, \bar t>4 r^2.$ Let $\Omega_1:= B_{2r}(\bar
x)\times [\bar t-2r^2, \bar t], \Omega:= B_r(\bar x)\times [\bar
t-r^2, \bar t].$ Suppose that there exist $\g\in \left[0, \frac 23
\right)$ and $\bar C>0$ such that
\begin{equation}\label{u6}
\sup_{\Omega_1} \e \xi^\e \leq \bar C \e^\g\quad \textrm{for
any}\;\; 0<\e<1\,.
\end{equation}
Then, for any $\s\in \left(\g, \frac 23(\g+1)\right)$, there
exists $C>0$ such that
\begin{equation}\label{u7}
\sup_{\Omega} \e\xi^\e \leq  C \e^\s\quad \textrm{for any}\;\;
0<\e<1\,.
\end{equation}
\end{lemma}

\noindent{\it Proof\,.\,\,} Fix any $\s\in \left(\g, \frac 23
(\g+1)\right)$. Define
\[G(u^\e):=\e^\s[2H_0 - (u^\e)^2]\,.\]
So, for some $\check C>0$, for any $0<\e<1$, $\phi_G^\e\leq \check
C \e^\g\quad \textrm{in}\;\; \Omega_1\,.$ We shall prove that, for
$\e\to 0^+$,
\begin{equation}\label{u18}
\sup_{\Omega}\phi_G^\e < \e^\s\,.
\end{equation} Note that from \eqref{u18} it follows that, for some
$C>0$,
\[\sup_{\Omega}\phi_G^\e \leq C \e^\s\quad \textrm{for any}\;\,0<\e<1,\]
so, \eqref{u7} follows (possibly taking a bigger $C>0$), and in
turn this inequality directly yields \eqref{u7}, by definition of $\phi^\e_G$.

\smallskip

Suppose, by contradiction, that \eqref{u18} is false and that
there exists a sequence $\{\e_n\}\subset (0,1)$ such that $\e_n\to
0$ as $n\to \infty$ and $\sup_{\Omega}\phi_G^{\e_n}\geq
\e_n^\s\,.$

\smallskip

Choose $\varphi\in C^\infty_0(\Omega_1)$ such that $0\leq
\varphi\leq 1, \varphi\equiv 1 $ in $\Omega$. Set
\[\tilde \phi^{\e_n}:= \phi_G^{\e_n} +\check C \e_n^\g \varphi.\]
Clearly, in $\big[B_{2r}(\bar x)\times\{\bar
t-2r^2\}\big]\cup\big[\pa B_{2r}(\bar x)\times (\bar t-2r^2, \bar
t] \big]$
\[\tilde \phi^{\e_n}\leq \sup_{\Omega_1}\phi_G^{\e_n}\leq \check C \e^\g\,.\]
Moreover,
\[\sup_{\Omega_1}\tilde \phi^{\e_n}\geq \sup_{\Omega}\tilde \phi^{\e_n} = \check C\e_n^\g+\sup_{\Omega}\phi_G^{\e_n}
\]
\[\geq \check C \e_n^\g +\e_n^\s>\check C \e_n^\g .\]
Therefore, $\max_{\bar \Omega_1}\tilde \phi^{\e_n}=\tilde
\phi^{\e_n}(x_n, t_n)$ for some $(x_n, t_n)\in \Omega_1\,.$ We
have:
\[\tilde \phi^{\e_n}(x_n, t_n)\geq \check C\e_n^\g +\e_n^\s>0,\]
so
\[\frac{\e_n^2}2|\nabla u^{\e_n}(x_n, t_n)|^2 \geq F(u^{\e_n}(x_n, t_n))+ G(u^{\e_n}(x_n, t_n))\]
\[+ \check C\e_n^\g(1- \varphi(x_n, t_n))+\e_n^\s.\]
Hence, for any $n\in \ene$,
\begin{equation}\label{u72}
 \frac{\e_n^2}2|\nabla
u^{\e_n}(x_n, t_n)|\geq \e_n^{1+\frac{\s}2}>0\,.
\end{equation}
Moreover, $\nabla \tilde \phi^{\e_n}(x_n, t_n)=0,$ thus, for some
constant $C_\varphi>0$, for any $n\in\ene$
\begin{equation}\label{u73}
|\nabla \phi_G^{\e_n}(x_n, t_n)|\leq C_{\varphi}\e_n^\g\,.
\end{equation}
We also have:
\[0\geq \Delta \tilde \phi^{\e_n}(x_n, t_n)-\pa_t \tilde \phi^{\e_n}(x_n, t_n)=\check C\e_n^\g[\Delta \varphi(x_n, t_n)-\pa_t\varphi(x_n, t_n)]\]
\[+\Delta \phi_G^{\e_n}(x_n, t_n)-\pa_t \phi_G^{\e_n}(x_n, t_n),\]
thus, for any $n\in\ene$,
\begin{equation}\label{u74}
\Delta \phi_G^{\e_n}(x_n, t_n)-\pa_t \phi_G^{\e_n}(x_n, t_n)\leq
C_{\varphi}\e_n^\g\,.\end{equation}

\smallskip
Let $\mathcal A^\e, \mathcal B^\e$ be defined as in Lemma
\ref{lemmau2}. We can find $\bar \e=\bar \e(\l)>0$ such that for
any $0<\e<\bar \e$
\[\frac{G''}{\e^2}-2\l <0\,.\]
For any $0<\e<\bar \e$, using \eqref{u72}-\eqref{u74}, we have
\[\mathcal A^\e \leq \tilde C\big[ \e_n^\g + \e_n^\g \e_n^{-1-\frac{\s}2}\big(|F'(u^{\e_n}(x_n, t_n))|+|G'(u^{\e_n}(x_n, t_n))|\big)\big]\,;\]
here and hereafter we always denote by $\tilde C$ possibly
different constants, independent of $n$ and $\e$. On the other
hand,
\[
\left(2\l-\frac{G''}{\e^2}\right)(F+G)\geq 0.\] Therefore, if
$|u^{\e_n}(x_n, t_n)|\geq \frac 12,$ then
\begin{equation}\label{u8}
\mathcal B^\e \geq \tilde C
\big[\frac{\e^{2\s}}{\e^2}+\frac{\e^\s}{\e^2}\big(|F'(u^{\e_n}(x_n,
t_n))|+|G'(u^{\e_n}(x_n, t_n))|\big)\big].
\end{equation}
If $|u^{\e_n}(x_n, t_n)|\leq \frac 1 2$,
\begin{equation}\label{u9}
\mathcal B^\e \geq \tilde C \frac{\e^\s}{\e^2}\big(\min_{|s|\leq
1/2}F(s) \big)\,.
\end{equation}
If $1\leq |u^{\e_n}(x_n, t_n)|\leq 1+ C\e^{\s_0}$, then
\begin{equation}\label{u10}
\mathcal B^\e \geq \tilde C \frac{\b_1 \e^{2\s}-\b_2
\e^{\s+\s_0}}{\e^2},
\end{equation}
for some $\b_1=\b_1(H_0)>0, \b_2=\b_2(F'')>0,$ for $\s_0>\s$
fixed.

\smallskip

Clearly, at least one inequality among \eqref{u8}, \eqref{u9},
\eqref{u10} holds for infinitely many $n\in \ene\,.$

Since $\g< \s < \frac 2 3(\g+1)$, we have $\s<\g+1-\frac{\s}2.$
So, for $n\in \ene$ large enough,
\[
\frac{\e_n^\s}{\e_n^2}\geq \tilde C \e_n^\g \e_n^{-1-\frac{\s}2}.
\]
Therefore, when $|u^{\e_n}(x_n, t_n)|\geq \frac 12$, for
$n\in\ene$ large enough, we have $\frac{\e_n^\s}{\e_n^2}\geq
\tilde C \e_n^\g,$ so
\begin{equation}\label{u75}
\e_n^{-2-\g+2\s}\leq \tilde C\,.
\end{equation}
However, $\g<\s<\frac 2 3(\g+1)<2.$ Therefore,
$-2-\g+2\s<-2-\g+\frac 43 (\g+1)=-\frac 2 3+\frac{\g}3<0.$ Hence
\eqref{u75} is impossible, for $n\in\ene$ large enough.

\smallskip

When $|u^{\e_n}(x_n, t_n)|\leq \frac 1 2$,
\[2 \min_{|s|\leq \frac 1 2}\{F(s)\}\frac{\e^\s}{\e^2} \leq \mathcal B^\e\leq  \mathcal A^\e \leq \tilde C (\e_n^\s + \e_n^\g \e_n^{-1-\frac{\s}2}).\]
Since $\g +1-\frac{\s}2>\s$, this yields
\[\frac{\e_n^{2\s}}{\e_n^2}\leq \tilde C \e_n^\g,\]
which is again impossible.

When $1\leq |u^{\e_n}(x_n, t_n)|\leq 1+ C\e^{\s_0}$, since $\s<\s_0$,
from \eqref{u10} we have
\begin{equation}\label{u67}
\mathcal B^\e \geq \tilde C \frac{\e_n^{2\s}}{\e^2}\,.
\end{equation}
Furthermore,
\begin{equation}\label{u68}
\mathcal A^\e\leq C(\e_n^{\g} +\e_n^\g\e_n^{-1-\frac{\s}2 +\s})\,.
\end{equation}
As above it is easily seen that \eqref{u67} and \eqref{u68} are in contrast. This completes the proof.
\hfill $\square$

\bigskip

Finally we are ready to prove Proposition \ref{propu1}.

\noindent{\it Proof of Proposition \ref{propu1}\,.\,\,} Let
$\Omega_0, \Omega$ be defined as in Lemma \ref{lemmau3}. Moreover,
set
\[\Omega_2:= B_{\frac 3 2 r}(\bar x)\times [\bar t- \frac 3 2 r^2, \bar t]\,.\]
By Proposition \ref{prop11b}, for some $C_0>0$,
\[\sup_{\Omega_1}\e \xi^\e \leq C_0 \e^0\quad \textrm{for all}\;\; 0<\e<1.\]
Thus, for any $0<\s<\frac 2 3$, for some $C_1>0$, by Lemma
\ref{lemmau3},
\[\sup_{\Omega_2}\e\xi^\e \leq C_1 \e^{\s}\quad \textrm{for all}\;\; 0<\e<1\,.\]
Hence, applying once more Lemma \ref{lemmau3}, for any
$0<\s<\frac{10}9$, for some $C>0$,
\begin{equation}\label{u14}
\sup_{\Omega}\e \xi^\e\leq C \e^\s\quad \textrm{for all}\;\;
0<\e<1\,.
\end{equation}
Now, the conclusion easily follows, choosing $1<\s<\frac{10}9$ in
\eqref{u14}. \hfill $\square$

\section{Asymptotic control of discrepancy for well-prepared initial
conditions}\setcounter{equation}{0}\label{wpic} In this Section we
prove an asymptotic control for the discrepancy $\xi^\e$, using
different methods from those in Section \ref{gic}. To this purpose
we need to assume that the initial conditions are properly well
prepared (see Subsection \ref{ic}) and the structure of the initial condition will emerge in the next two paragraphs.

\subsection{One-dimensional profile}\label{1d} Now we study problem
\eqref{e38}, where $\varphi:[0,1] \to(0,\infty)$ is an
increasing convex smooth function such that $\varphi(0)\ge
1$ and $\varphi^\prime(0)=0$. Let us define the energy
\[\mathcal E(h_\e):=\int_{0}^{1} \left(\frac 1 2 h'_\e\,^2 + F_{\e}(h_\e)\right)\varphi(\t)d\t\,.\]

By the same arguments as in the proof of Proposition 3.1 in
\cite{PisPons} it is possible to show next

\begin{lemma}\label{lemma2u}
For any $\e>0$ there exists a unique solution $h_\e$ to problem
\eqref{e38}. Furthermore, $h_\e$ is increasing and concave in
$[0,1]$, and there holds:
\begin{equation}\label{e39}
\mathcal E_\e(h_\e)\le \frac {C_1}{\e},\quad 0<h_\e'\le \frac
{C_1}{\e}\,,
\end{equation}
for some positive constant $C_1$ independent of $\e$.
\end{lemma}

The following lemma gives the main property of the profile function $h_\e$.
\begin{lemma}\label{lemma16} For any $\e>0$ let $h_\e$ be the unique solution to problem
\eqref{e38}, and still denote by $h_\e$ its odd reflection. Then
\eqref{e161} holds true.
\end{lemma}

\noindent{\it Proof\,\,.} Clearly we may assume $\t\geq 0$. From
\eqref{e38} we get
\begin{equation}\label{e181}
\frac{d}{d\t}\left(\frac 1 2h_\e'\,^2-F_\e(h_\e)\right)=h_\e'
\big[h_\e''-F_\e'(h_\e)\big]=-\frac{\varphi'}{\varphi}h_\e'\,^2\leq
0\quad \textrm{in}\;\;(0,1)\,.
\end{equation}
Still denote by $\varphi$ its even reflection. Since $h_\e$ is odd
and $\varphi$ is even, $\varphi'(\t)>0$ for all $\t>0$, from
\eqref{e181} we get
\begin{equation}\label{e180}
\e\left(\frac 1 2[h_\e'(\t)]^2- F_\e(h_\e(\t))\right)\le\e\frac 12
[h_\e'(1)]^2+\e\int_{|\t|}^{1}\frac{\varphi'}{\varphi}h_\e'\,^2
ds\quad \textrm{in}\;\;(-1,1)\,.
\end{equation}
Furthermore, using \eqref{e39} we have
\begin{equation}\label{e134}
\begin{split}
\e \int_0^{1} \frac{\varphi'}{\varphi}(h')^2 ds\le
\left(\e\int_0^{1}\frac{(h')^2}{\varphi^2}ds\right)^{\frac 1 2}
\left(\e\int_0^{1}(h')^2 \varphi'^2 ds \right)^{\frac 1 2}\\
\le C_1 \left(\e\int_0^{1}(\varphi')^2(h')^2 ds\right)^{\frac 1
2}.\hspace{1.5 cm}
\end{split}
\end{equation}
We shall prove that
\begin{equation}\label{e135}
|h'_\e|\le \frac{o(1)}{\sqrt \e}\quad \textrm{if}\;\; \t>\sqrt
\e\,.
\end{equation}
To this purpose, note that since $h_\e$ is concave in $(0,1)$, we
have:
\begin{equation}\label{e136}
0<h'_\e(\t)\le \frac{1-h_\e(k\e)}{\sqrt \e-k\e}\quad \textrm{for
all}\;\; \t\in[\sqrt\e, 1]\,,
\end{equation}
for each $k\in (0,\infty)$ and for $\e\in (0,1)$ so small that
$k\e\leq 1$. The function $v_{\e}(s):= h_{\e}(s \e)\,\, (s\in [0,
\frac 1{\epsilon}])$ solves
\[v_\e '' + \e \frac{\varphi'(s\e)}{\varphi(s\e)}v'_\e = f(v),\]
and it is easy to see that it converges as $\e\to 0^+$ in
$C^2_{loc}(\re)$ to a monotone increasing solution $v=v(s)$ of
equation
\[v''-f(v)=0\quad \textrm{in}\;\;(-\infty,\infty)\,.\]

Now, note that
\[\e^2\left\{\frac 1 2h_\e'(\t)\,^2-
F_\e(h_\e(\t))\right\}=\e^2 \int_\t^1
\frac{\varphi'(s)}{\varphi(s)}h_\e'^2(s) ds +\e^2\left\{\frac 1
2h_\e'(1)\,^2- F_\e(h_\e(1))\right\}\,.\] In view of concavity of
$h_\e$, it is easily seen that
$$\e^2\left\{\frac 1
2h_\e'(1)\,^2- F_\e(h_\e(1))\right\}\to 0\quad\textrm{as}\;\;
\e\to 0^+,
$$
while
\[\e^2 \left| \int_\t^1
\frac{\varphi'(s)}{\varphi(s)}h_\e'^2(s) ds \right|\le C \e^2
\mathcal E(h_\e)\le C C_1\e,
\]
in view of \eqref{e39}. Hence, from \eqref{e181} we have
$$\e^2\left\{\frac 1 2h_\e'(\t)\,^2- F_\e(h_\e(\t))\right\}\to 0
\quad\textrm{as}\;\; \e\to 0^+\quad \textrm{for any}\;\; \t\in
[-1,1]\;\;.$$ As a consequence,
$$v_\e'(0)=\e h_\e'(0)\to F(0)>0\quad \textrm{as}\;\; \e\to 0^+.$$
Thus, $v'(0)=F(0)>0, v$ is bounded and strictly increasing, hence
$v(s)\to 1$ as $|s|\to \infty$. Since $h_\e(k\e)=v_\e(k)\to v(k)$
as $\e\to 0,$ this combined with \eqref{e136} gives \eqref{e135}.

\smallskip

Observe that, since $\varphi$ is smooth and $\varphi^\prime(0)=0$, there exists
$C>0$ such that
\begin{equation}\label{e137}
\varphi'(\t)\le C\t \quad \textrm{in}\;\; (0,1)\,.
\end{equation}

\smallskip

Inequalities \eqref{e134}, \eqref{e135} and \eqref{e137} yield
\begin{equation}\label{e138}
\begin{split}
\e \int_0^{1}(\varphi')^2 (h_\e')^2 ds \le C\e\Big(
\int_0^{\sqrt\e} \frac{\t^2}{\e^2}ds \hspace{.5 cm}\\+
\int_{\sqrt\e}^{1}\frac{o(1)}{\e}ds \Big)\le C[\sqrt\e + o(1)]\to
0\quad \textrm{as}\;\; \e\to 0^+\,.
\end{split}
\end{equation}
From \eqref{e180}, \eqref{e134}, \eqref{e135} and \eqref{e138} the
conclusion follows. \hfill $\square$

\subsection{Well prepared initial conditions}\label{ic}
Now we are ready to define well prepared initial conditions
$u_0^\e$. To be specific, we assume that the initial conditions
$u_0^\e$ are in the form \eqref{e71}. In addition we assume that
$u_0^\e$ and the corresponding $z^\e_0\equiv z^\e(\cdot,0),
\mu^\e_0\equiv \mu^\e(\cdot,0)$ given by \eqref{e71},\eqref{e8i}
satisfy
\[\left\{
\begin{array}{l}
(i)\;\, (H_1)\,\, \textrm{is satisfied}\,;\\
(ii)\, z^\e(\cdot,0) \in C^2(M) \quad \hbox{and} \quad |\nabla z^\e(x,0)|\le 1\;\,\hbox{for all}\;\; x\in M\,. \vspace{.2 cm}\\
\end{array}
\right. \leqno(H_1^*)
\]
The construction of such a $u_0^\e$ is quite standard
(see \cite{Modica}, \cite{Ilm1}). %Let $\check d$ be a smoothing of $\tilde d(x, \Sigma_0)$ that
%equals $\tilde d(x, \Sigma_0)$ near $\Sigma_0$, and which
%satisfies $|\nabla \check d|\le 1$.
Since we assume $\Sigma_0$ to be smooth (at least $C^3$), there is a small tubular neighboorod $U_\delta\supset \Sigma_0$ of size $4\delta>0$ such that the distance function $d(x,\Sigma_0)$ is smooth in $U_\delta$ (at least $C^2$). Let now $\Psi \in C^\infty(\mathbb{R})$ an odd increasing function such that $\Psi(s)=s$ whenever $|s|<\delta$, $|\Psi(s)|=2\delta$ for $|s| \geq 4 \delta$ and such that $\Psi^{\prime\prime}\leq 0$ for $s>0$. Then it is direct to see that if $\tilde{d}(x, \Sigma_0)$ is the signed distance from $\Sigma_0$, then $z^\e(x,0)=\Psi(\tilde{d}(x), \Sigma_0)$ is a globally smooth function (as smooth as the distance is near $\Sigma_0$), is constant far from $\Sigma_0$ and the corresponding $u_0^\e$ given by

\[ u_0^\e(x)=h_\e( \Psi(\tilde{d}(x, \Sigma_0))) \, , \]
 with $h^\e$ as in Lemma \ref{lemma2u}, satisfy assumption $(H_1^*)$ above.

%Concerning the initial conditions $u^\e_0$, well prepared data
%will be now of the form  (see Subsection \ref{ic}):
%\begin{equation}
%\label{wellprepdata}
%u_0^{\e}(x):= h^\e\left( \Psi ( \tilde d(x, \Sigma_0) )
%\right)\, \quad (x\in M)\,,
%\end{equation}
%where $E_0\subset M$ is an open bounded subset with smooth
%boundary $\Sigma_0:=\pa E_0$. Now, $d(x, \Sigma_0)$  is the
%Riemannian distance of $x\in M$ to $\Sigma_0$, while the signed
%distance $\tilde d(x,\Sigma_0)$ is defined accordingly and $\Psi(\cdot)$ is a suitable smoothed and $1-$Lipschitz truncation of the identity, which makes $u_0^\e$ constant far from $\Sigma_0$, where the distance function is possibly singular.

 Let $\mathcal F(s):=\int_0^s \sqrt{F(\t)}d\t$. We have
that $u_0^\e\to u_0^0$, and $\mathcal F(u_0^\e)\to \mathcal F
(u_0^0)$ as $\e \to 0^+$, uniformly in $M\setminus \Sigma_0$ and
in $L^1_{loc}(M)$; furthermore, for $c=\int_{-1}^1 \sqrt{F(s)}ds \, ,$
\[\int_{\Omega}\left(\frac \e 2|\nabla
u_0^\e|^2+\frac{F(u_0^\e)}{\e}\right) d \mathcal V\to c \, |\nabla
u_0^0|(\Omega)\quad \textrm{as}\;\,\e \to 0^+,\;\; \textrm{for
every}\;\,\Omega\subset\subset M\,,
\]
where $|\nabla u_0^0|(\Omega)$ denotes the total variation of
$u_0^0$ in $\Omega$, whence $\mathcal F(u_0^\e)\to \mathcal F
(u_0^0)$ as $\e \to 0^+$ weakly-* in $BV_{loc}(M)$.

Note that if $(H_0), (H_1^*)$ are satisfied, then the unique
solution to problem \eqref{e1}-\eqref{e1a}, which exists by
Proposition \ref{prop12}, also verifies
\begin{equation}\label{e121}
-1 < u^\e < 1\quad\textrm{in}\;\; M\times (0,\infty)\,.
\end{equation}
This follows by maximum principle, since $|u_0^\e|\leq 1$.
In addition, since we assume $u_0^\e \in C^2(M)$, parabolic regularity theory also implies
$\nabla u^\e \in C^0(M \times [0,\infty))$.

\subsection{Asymptotic control of discrepancy}\label{discrep}
In order to show \eqref{u13} we need to prove preliminarily that
\eqref{e52} follows from $(H_0), (H_1^*)$.

\begin{lemma}\label{lemma6}
Let assumptions $(H_0), (H_1^*)$ be satisfied. Let $u^\e$ be a
solution to equation \eqref{e1}; suppose that \eqref{e55b} and
\eqref{e51} hold true. Then inequality \eqref{e52} is satisfied.
\end{lemma}

\noindent{\it Proof\,.\,\,} Define $w^\e:= |\nabla z^\e|^2\,,$ and note that, as already proved above, $w^\e \in C^0 (M\times[0,\infty)) \cap C^\infty(M\times(0,\infty))$.
From \eqref{e1}, \eqref{e53}, \eqref{e53a} and \eqref{e38} we
deduce that
\begin{equation}\label{e54}
\begin{split}
\pa_t w^\e = \Delta w^\e - 2 |\Hess z^\e|^2-\l w^\e \hspace{1.5 cm}\\
-\frac{f(u^\e)}{\e^2 h'_\e(z^\e)}\langle\nabla z^\e, \nabla w^\e
\rangle
-\frac{w^\e}{\e^2}(w^\e-1)\left[f'(u^\e)-f(u^\e)\frac{h''(z^\e)}{(h'(z^\e))^2}
\right]\\-\frac{\varphi'}{\varphi}(z^\e)\langle\nabla z^\e, \nabla
w^\e \rangle-\left(\frac{\varphi'}{\varphi}\right)'(z^\e)w^\e
\quad \textrm{in}\;\;M\times (0,\infty)\,. \hspace{.7 cm}
\end{split}
\end{equation}
Note that, in view of \eqref{e51}, the function $\overline w\equiv
1$ is a supersolution to equation \eqref{e54}. Note that in view
of \eqref{e53}, from results in \cite{Dod} comparison principles
can be easily obtained. Hence, from \eqref{e55b} and comparison
principles inequality \eqref{e52} follows. \hfill $\square$

Now we can prove the following proposition.
\begin{proposition}\label{prop2}
Let assumptions of Lemma \ref{lemma6} be satisfied. Then
\begin{equation}\label{u13b}
\limsup_{\e\to 0^+} \sup_{(x,t)\in M\times(0,\infty)}
\xi^\e_t(x)\leq 0\,.
\end{equation}
Hence, in particular, \eqref{u13} holds true.
\end{proposition}

\noindent {\it Proof\,.\,\,} From \eqref{e181} and \eqref{e52}, we
get:
\begin{equation}\label{e56}
\begin{split}
\frac 1 2|\nabla u^\e|^2- F_\e (u^\e)=\frac 1 2
[h_\e'(z^\e)]^2|\nabla z^\e|^2-F_\e(h_\e(z^\e)) \\
\le \frac 1 2 [h'_\e(z^\e)]^2-
F_\epsilon(h_\epsilon(z_\epsilon))\quad \textrm{in}\;\;M\times (0,
\infty)\,.\hspace{.3 cm}
\end{split}
\end{equation}
From Lemma \ref{lemma16} the conclusion follows. \hfill $\square$

\section{Uniform energy bounds}
\setcounter{equation}{0} This section is devoted to proof the {\it
local almost monotonicity formula} \eqref{e60} and to derive from
apriori estimates some compactness properties of the family of
solutions $u^\e$ as $\e \to 0$ both in $BV_{loc}$ and in
$L^1_{loc}$.

\subsection{Local almost monotonicity formula}\label{mf}
The argument to prove \eqref{e60} is a modification of the one in \cite{Ilm1}, localizing the  estimate at suitably small scale so to reabsorbe the perturbation terms coming from the curved background. This, combined with the locally uniform control of the positive part of the discrepancy from the previous sections, will give the final result.

At first, recall the next lemma (see Lemma
6.6 in \cite{Ilm2}).
\begin{lemma}\label{lemma14}
Let $\varphi\in C^2_c(M;[0,\infty))$. Then
\[\frac{|\nabla \varphi|^2}{\varphi}\le 2 \max_{\{\varphi>0\}}|\Hess \varphi|\quad \textrm{in}\;\;\{\varphi>0\}\,.\]
\end{lemma}

Then, the next lemma will give a a sort
of {\it Brakke's} inequality for $d\mu^\e_t$.

\begin{lemma}\label{lemma1}
Let $u^\e$ be a solution to equation \eqref{e1}. Let $\phi\in
C^{2,1}_{x,t}(M\times (0,\infty);\re_+)$ with $\suppo
\phi(\cdot\,, t)$ compact for every $t\in (0,\infty)$. Then
\begin{equation}\label{e4}
\begin{split}
\frac{d}{dt} \int_{M} \phi E^\e d\mathcal V(x) dt  \hspace{2.5 cm}\\
=\int_{M}\Big\{(\pa_t\phi - \Delta \phi)E^\e + (\Hess \phi)(\nabla
u^\e,\nabla u^\e) - \phi(\pa_t u^\e)^2\Big\}d\mathcal V(x)\,,
\end{split}
\end{equation}
and
\begin{equation}\label{e9}
\begin{split}
\frac{d}{dt} \int_{M} \phi E^\e d \mathcal V(x) dt  \hspace{2.8 cm}\\
=\int_{M}\Big\{(\pa_t\phi + \Delta \phi)E^\e - (\Hess \phi)(\nabla
u^\e,\nabla u^\e)\Big\}d\mathcal V(x)  \hspace{1 cm}\\+ \int_{M}
\frac{\langle \nabla \phi, \nabla u^\e\rangle^2}{\phi}d\mathcal
V(x) - \int_{M} \phi\left(\pa_t u^\e +\frac{\langle \nabla \phi,
\nabla u^\e\rangle}{\phi} \right)^2d\mathcal V(x)\,.
\end{split}
\end{equation}
for all $t>0$.
\end{lemma}

Note that the last two integrals in equality \eqref{e9} are well
defined, in view of Lemma \ref{lemma14}.

\medskip

\noindent{\it Proof\,.\,\,} By \eqref{e1}, in $M\times
(0,\infty)$,
\begin{equation}\label{e6}
\phi\pa_t E^\e = \phi\langle \nabla \pa_t u^\e , \nabla
u^\e\rangle + \phi\pa_t u^\e\Delta u^\e - \phi(\pa_t u^\e)^2\,,
\end{equation}
and
\begin{equation}\label{e7}
\langle \nabla \phi,\nabla E^\e\rangle=\frac 1 2\langle\nabla
\phi, \nabla \langle \nabla u^\e, \nabla u^\e\rangle \rangle-
\pa_t u^\e\langle \nabla \phi , \nabla u^\e\rangle+ \langle\nabla
u^\e, \nabla \phi\rangle\Delta u^\e\,.
\end{equation}
In view of \eqref{e6}-\eqref{e7} we get
\begin{equation}\label{e10}
\begin{split}
\phi \pa_t E^\e -\langle \nabla \phi,\nabla E^\e\rangle=
\phi\langle \nabla \pa_t u^\e , \nabla u^\e\rangle + \phi\pa_t
u^\e\Delta u^\e - \phi(\pa_t u^\e)^2\\
-\frac 1 2\langle\nabla \phi, \nabla \langle \nabla u^\e, \nabla
u^\e\rangle \rangle+ \pa_t u^\e\langle \nabla \phi , \nabla
u^\e\rangle-\langle\nabla u^\e, \nabla \phi\rangle\Delta u^\e\,.
\end{split}
\end{equation}
From \eqref{e10}, integrating by parts and using \eqref{e5}, it
follows:
\[
\frac{d}{dt}\int_{M} \phi E^\e d\mathcal V(x) = \int_{M} \pa_t\phi
E^\e d\mathcal V(x)\]\[ + \int_{M}\Big\{\langle\nabla \phi, \nabla
E^\e \rangle -\phi (\pa_t u^\e)^2 + \phi\pa_t u^\e\Delta u^\e +
\phi \langle \nabla\pa_t u^\e, \nabla u^\e \rangle \]\[-\frac 1
2\langle \nabla\phi, \nabla\langle \nabla u^\e, \nabla u^\e\rangle
\rangle + \pa_t u^\e \langle\nabla \phi, \nabla u^\e \rangle -
\Delta u^\e\langle \nabla\phi, \nabla u^\e\rangle\Big\}d \mathcal
V(x) \]\[= \int_{M}(\pa_t\phi - \Delta\phi) E^\e d \mathcal V(x)
-\int_{M} \phi(\pa_t u^\e)^2 d \mathcal V(x) \]\[-
\int_{M}\Big\{\frac 1 2\langle \nabla \phi, \nabla \langle\nabla
u^\e, \nabla u^\e\rangle \rangle - \langle \nabla u^\e, \nabla
\langle \nabla \phi, \nabla u^\e\rangle\rangle\Big\} d \mathcal
V(x)\]\[ =\int_{M}(\pa_t\phi - \Delta\phi) E^\e d \mathcal V(x)
-\int_{M} \phi(\pa_t u^\e)^2 d \mathcal V(x) \]\[-\int_{M}\Big\{
\frac 1 2\Hess (u^\e)\langle \nabla u^\e, \nabla\phi\rangle +
\frac 1 2\Hess (u^\e)\langle\nabla u^\e, \nabla \phi\rangle \]\[-
\Hess (u^\e)(\nabla u^\e, \nabla\phi)\,+\,\Hess (\phi)(\nabla
u^\e, \nabla u^\e)\Big\}d \mathcal V(x)\,\]\[= \int_{M}(\pa_t\phi
- \Delta\phi) E^\e d \mathcal V(x) -\int_{M} \phi(\pa_t u^\e)^2 d
\mathcal V(x) \]\[+\int_{M} \Hess (\phi)(\nabla u^\e,\nabla u^\e)
d \mathcal V(x)\,.
\]
Hence \eqref{e4} has been verified. Equality \eqref{e9} can be
shown analogously, using
\[
\phi \pa_t E^\e + \langle \nabla \phi,\nabla E^\e\rangle=
\phi\langle \nabla \pa_t u^\e , \nabla u^\e\rangle + \phi\pa_t
u^\e\Delta u^\e - \phi(\pa_t u^\e)^2\]\[ +\frac 1 2\langle\nabla
\phi, \nabla \langle \nabla u^\e, \nabla u^\e\rangle \rangle -
\pa_t u^\e\langle \nabla \phi , \nabla u^\e\rangle+\langle\nabla
u^\e, \nabla \phi\rangle\Delta u^\e\,,
\]
instead of \eqref{e10}.\hfill $\square$

\bigskip
Next we are going to chose a precise test function in the formulas obtained above. The following lemma gives auxiliary identities which will be useful in this direction.
\medskip

\begin{lemma}\label{lemma20}
Let $K\subset M$ be a compact subset, $y\in K, s>0$. Let $\hat \zeta\in C^2([0,\infty))$ such that
\begin{equation}\label{e11}
|\hat \zeta|\le 1\,,\;\; |\hat \zeta'|\le 1\,,\, |\hat \zeta''|\le 1\quad
\textrm{in}\;\; [0,\infty)\,,
\end{equation}
\begin{equation}\label{e12}
\hat \zeta = \left\{
\begin{array}{ll}
1 & \textrm{in} \;\;  [0,R_0^2/4)
\\& \\
0  & \textrm{in} \;\; [R_0^2 , \infty)\,,
\end{array}
\right.
\end{equation}
where $R_0:= \frac 1 2\inf_{y\in K}\inj(y)\,.$ Define
\[\hat\eta(\rho,t):=[(s-t)]^{-\frac{N-1}2} e^{-\frac{\rho}{4(s-t)}}\quad (\rho\geq 0, 0\le t<s)\,.\]
Let
\begin{equation}\label{u52}
\eta(x,t):= \hat \eta\big(d^2(x), t\big)\quad (x\in M, 0\leq
t<s)\,.
\end{equation}
\[\zeta(x):= \hat \zeta (d^2(x))\quad (x\in M)\,.\]
Then, for all $x\in M, 0\leq t<s,$
\begin{equation}\label{u53}
\nabla (\eta\zeta)(x,t)= \pa_\rho (\hat \eta \hat \zeta)(d^2(x), t)\nabla d^2(x),
\end{equation}
\begin{equation}\label{u54}
\begin{split}
\Hess [(\eta\zeta)](d^2(x),t)(X,X) = \pa_{\rho\rho} (\hat \eta
\hat\zeta)(d^2(x), t)\big|\langle \nabla d^2(x)  , X  \rangle
\big|^2 \\+ \pa_\rho (\hat \eta \hat\zeta)(d^2(x),
t)\Hess[d^2(x)](X,X)\quad (X\in T_xM)\,;\hspace{.5 cm}
\end{split}
\end{equation}
\begin{equation}\label{u55}
\Delta (\eta\zeta)(x,t)= \pa_{\rho\rho} (\hat\eta\hat
\zeta)(d^2(x), t)|\nabla d^2(x)|^2 + (\pa_\rho
\hat\eta\hat\zeta)(d^2(x), t)\Delta d^2(x)\,;
\end{equation}
\begin{equation}\label{u56}
\begin{split}
(\pa_t +\Delta)(\eta\zeta)(x,t)\hspace{3.5 cm}\\ = \frac{(\eta\zeta)(x,t)}{s-t}\left[\frac{N-1}2-\frac{d^2(x)}{4(s-t)}+\frac{|\nabla d^2(x)|^2}{16(s-t)}
-\frac 1 4 \Delta d^2(x)\right] \hspace{1 cm}\\
+ (2\pa_{\rho} \hat\eta \pa_\rho \hat\zeta +\hat \eta
\pa_{\rho\rho}\hat\zeta)(d^2(x), t) |\nabla d^2(x)|^2 + (\hat
\eta\pa_\rho\hat \zeta) (d^2(x), t)\Delta d^2(x)\,.\hspace{.3 cm}
\end{split}
\end{equation}
\end{lemma}

\begin{remark}\label{ossu1} It is straightforward to check that
\begin{equation}
\pa_\rho \hat \eta=-\frac{\hat \eta}{4(s-t)},  \;\;
\pa_{\rho\rho}\frac{\hat \eta}{16(s-t)^2}\,,\;\; (\pa_\rho \hat
\eta)^2 - \hat \eta \pa_{\rho\rho}\hat \eta \,=\,0\,.
\end{equation}
\end{remark}

\noindent{\it Proof of Lemma \ref{lemma20}\,.} For any $x\in M, 0\leq t<s$ we have:
\[\nabla (\eta\zeta)(x,t)=\pa_\rho(\hat\eta\hat\zeta)(d^2(x), t)\nabla d^2(x)\,.\]
Hence, for any $X\in T_xM$,
\[\Hess[(\eta\zeta)](x,t)(X,X)=\langle \cd_X \nabla (\eta\zeta), X\rangle\]\[= \langle \cd_X (\pa_\rho \hat\eta\hat\zeta)(d^2(x),t) \nabla d^2(x), X \rangle
=\pa_{\rho\rho}(\hat\eta \hat \zeta)(d^2(x), t)|\langle X ,\nabla
d^2(x)\rangle|^2\]\[+  \pa_\rho(\hat \eta \hat \zeta)(d^2(x),
t)\Hess[d^2(x)](X,X)\,.\] Passing to the trace we get
\[\Delta(\eta\zeta)(x,t)=\pa_{\rho\rho}(\hat\eta\hat\zeta)(d^2(x), t)|\nabla d^2(x)|^2+\pa_\rho(\hat\eta\hat\zeta)(d^2(x), t)\Delta d^2(x)\,.\]

\smallskip
Furthermore,
\begin{equation}\label{u57}
\pa_t \hat\eta =
-\frac{1-N}{2}\frac{\hat\eta}{s-t}+\hat\eta\left(-\frac{\rho}{4(s-t)^2}\right)=\frac{\hat\eta}{s-t}\left[\frac{N-1}2-\frac{\rho}{4(s-t)}\right]\,.
\end{equation}
Note that
\begin{equation}\label{u58}
(\pa_t +\Delta)(\eta\zeta)(x,t)=\zeta(\pa_t +\Delta)\eta(x,t) + 2 \nabla\eta(x,t)\nabla \zeta(x,t)+\eta\Delta\zeta(x)\,;
\end{equation}
moreover,
\begin{equation}\label{u59}
\nabla\eta(x,t)=\pa_\rho\hat\eta(d^2(x), t)\nabla d^2(x)\,,
\end{equation}
and
\begin{equation}\label{u60}
\Delta\eta(x,t)=\pa_{\rho\rho}\hat\eta(d^2(x), t)|\nabla d^2(x)|^2 + \pa_\rho\hat\eta(d^2(x),t)\Delta d^2(x)\,.
\end{equation}
By \eqref{u58}-\eqref{u60},
\[(\pa_t +\Delta)\eta(x,t)=\frac{\hat\eta}{s-t}\left[\frac{N-1}2-\frac{d^2(x)}{4(s-t)}+ \frac{|\nabla d^2(x)|^2}{16(s-t)}-\frac 1 4 \Delta d^2(x) \right],  \]
so
\[(\pa_t +\Delta)(\eta\zeta)(x,t)\]\[=\frac{(\hat\eta\hat\zeta)(d^2(x), t)}{s-t}\left[-\frac 1 2 +\frac{d^2(x)}{4(s-t)}(-1+|\nabla d(x)|^2 )+\frac N 2-\frac 1 4\Delta d^2(x) \right]   \]
\[+ \big(2\pa_\rho\hat\eta \pa_\rho\hat\zeta +\hat\eta\pa_{\rho\rho}\hat\zeta\big)(d^2(x),t)|\nabla d^2(x)|^2+(\hat\eta\pa_\rho\hat\zeta)(d^2(x), t)\Delta d^2(x)\,. \]
This implies \eqref{u56}, since
\[
|\nabla d(x)|^2 \leq 1\,.
\]
This completes the proof. \hfill $\square$

\bigskip

In order to proceed we start rewriting some terms in \eqref{e9}. Note that
\[ \frac{|\langle \nabla\phi, \nabla u^\e \rangle  |^2}{\phi}-\Hess \phi (\nabla u^\e, \nabla u^\e)\]\[= \frac{| \langle \nabla d^2(x)  ,
 \nabla u^\e \rangle  |^2}{\eta\zeta}
(\hat\zeta\pa_\rho \hat\eta +\hat\eta \pa_\rho\hat\zeta)^2-
(\hat\zeta \pa_\rho\hat\eta
+\hat\eta\pa_\rho\hat\zeta)\Hess[d^2(x)](\nabla u^\e, \nabla
u^\e)\]\[- |\langle \nabla d^2(x) , \nabla u^\e \rangle
|^2(\hat\zeta \pa_{\rho\rho}\hat\eta + 2 \pa_\rho\hat\zeta
\pa_\rho \hat\eta + \hat \eta \pa_{\rho\rho}\hat\zeta)\]
\[=\frac{\hat\eta \hat\zeta}{4(s-t)}\Hess[d^2(x)](\nabla u^\e,
\nabla u^\e)-\hat\eta \pa_\rho\hat\zeta\Hess[d^2(x)](\nabla u^\e,
\nabla u^\e)\]\[ + |\langle \nabla d^2(x), \nabla u^\e\rangle|^2
\Big[\frac{\hat\zeta^2(\pa_\rho\hat \eta)^2+ 2
\hat\eta\pa_\rho\hat\eta \hat\zeta\pa_\rho\hat\zeta+\hat\eta^2
(\pa_\rho\hat\zeta)^2}{\hat
\eta\hat\zeta}\]\[-\hat\zeta\pa_{\rho\rho}\hat\eta-2
\pa_\rho\hat\eta\pa_\rho\hat\zeta-\hat\eta\pa_{\rho\rho}\hat\zeta\Big]\,.
\] So, we can infer the following lemma.

\begin{lemma}\label{lemma3} Let $\eta$ and $\zeta$ be as in Lemma
\ref{lemma20}, and $\phi:=\eta\zeta$. Then
\[\int_{M}\left[\e\frac{| \langle \nabla \phi, \nabla u^\e\rangle |^2}{\phi} -\e \Hess \phi (\nabla u^\e, \nabla u^\e)\right]d\mathcal V(x)
=\e\int_M  \frac{\eta\zeta}{2(s-t)}|\nabla u^\e|d\mathcal V(x)
\]\[+ \e \int_M \frac{\eta\zeta}{2(s-t)}\left[\frac 1 2
\Hess[d^2(x)](\nabla u^\e, \nabla u^\e)-|\nabla u^\e|^2
\right]d\mathcal V(x)  \]
\[+\e\int_M |\langle \nabla d^2(x), \nabla u^\e \rangle   |^2 \eta \left(\frac{\pa_\rho \hat \zeta}{\hat\zeta}-\pa_{\rho\rho}\hat\zeta \right) d\mathcal V(x)
\]\[-\e \int_M \eta\pa_\rho\hat\zeta \Hess[d^2(x)](\nabla u^\e,
\nabla u^\e)d\mathcal V(x)\,.\]
\end{lemma}

\bigskip
\smallskip

From Lemmas \ref{lemma1}, \ref{lemma20}, \ref{lemma3}, neglecting one negative term we
immediately get
\begin{lemma}\label{propu2}
Let $\phi:= \eta \zeta$. Then
\[\frac{d}{dt}\int_M \phi\e E^\e d\mathcal V(x) \leq \int_M \frac{\phi}{2(s-t)}\big[\e|\nabla u^\e|^2-\e E^\e \big]d\mathcal V(x)   \]
\[+ \e \int_M \frac{\phi}{2(s-t)}\left[ \frac 1 2 \Hess[d^2(x)](\nabla u^\e, \nabla u^\e) -|\nabla u^\e|^2 \right]d\mathcal V(x)\]\[ + \e \int_M \frac{\phi}{2(s-t)}\left[N-\frac 1 2\Delta d^2(x) \right]
E^\e d\mathcal V(x) + \e \int_M E^\e \eta\pa_\rho \hat\zeta \Delta d^2(x)\]
\[+\e \int_M E^\e |\nabla d^2(x)|^2 (2\pa_\rho\hat\eta\pa_\rho\hat\zeta +\hat\eta\pa_{\rho\rho}\hat\zeta)d\mathcal V(x)
\]\[- \e \int_M \eta\pa_\rho\hat\zeta\Hess[d^2(x)](\nabla u^\e,
\nabla u^\e)d\mathcal V(x)   \]
\[+\e\int_M |\langle \nabla d^2(x)  , \nabla u^\e \rangle  |^2\eta\left(\frac{(\pa_\rho\hat\zeta)^2}{\hat\zeta}-\pa_{\rho\rho}\hat\zeta \right)d\mathcal V(x)\,.\]
\end{lemma}

Finally, we single out the discrepancy term and enstimate all the others to obtain the key result of this section.

\begin{proposition}\label{lemma5}
Let assumption $(H_0)$ be satisfied. Let $u^\e$ be the solution to
problem \eqref{e1}-\eqref{e1a}. Suppose that \eqref{e121b} and
\eqref{e170} with $\t=0$ hold true. Let $K\subset M$ be a compact
subset, $y\in K, s>0.$ Let $\phi:=\eta\zeta$ with $\eta$ and
$\zeta$ as in Lemma \ref{lemma20}. Then for every $\e\in (0,1)$
\begin{equation}\label{e47}
\begin{split}
\frac{d}{dt}\int_{M} \phi(x,t)d\mu^\e_t(x)\,\le\, \frac 1{2(s-t)}
\int_{M} \phi\e\big\{|\nabla u^\e|^2-E^\e \big\}d \mathcal V(x)
\\+\frac {C_3}{(s-t)^{1/2}}\int_{M} \phi
d\mu^\e_t(x)+ C_4\quad \textrm{for
all}\;\;0<t<s\,.\hspace{.3 cm}
\end{split}
\end{equation}
for some positive constants $C_3$ and $C_4$ depending only on $K$.
\end{proposition}

\noindent{\it Proof\,.\,\,} By Lemma \ref{lemmau1} and Lemma
\ref{propu2},
\begin{equation}\label{e32}
\begin{split}
\frac{d}{dt} \int_{M} \phi \e E^\e d \mathcal V(x)
\le \frac 1{2(s-t)} \int_M \phi d\xi^\e_t + C\int_M \phi \frac{d^2(x)}{s-t}\e E^\e d\mathcal V(x)\\
+ C_{R_0}\int_M \e E^\e \eta\left(|\pa_\rho
\hat\zeta|+\frac{|\pa_\rho\hat\zeta|}{\hat\zeta}+|\pa_{\rho\rho}\hat\zeta|\right)d\mathcal
V(x) \hspace{1 cm}\\
+C_{R_0}\int_M\e E^\e \eta\left(\frac{d^2(x)}{s-t}|\pa_\rho\hat\zeta|+|\pa_{\rho\rho}\hat\zeta|\right)d\mathcal V(x) \hspace{1.2 cm}\\
\leq \frac 1{2(s-t)}\int_M \phi d\xi^\e_t + C \int_M
\phi\frac{d^2(x)}{s-t}\e E^\e d\mathcal V(x) \hspace{1.2 cm} \\
+\tilde C_{R_0}\int_{\frac{R_0}2\leq d(x)\leq R_0}\e E^\e
d\mathcal V(x)\,,\hspace{1.5 cm}
\end{split}
\end{equation}
for some positive constants $C, C_{R_0}, \tilde C_{R_0}$ independent of $\e$.

Now, note that
\begin{equation}\label{u62}
\begin{split}
 \int_M \phi\frac{d^2(x)}{s-t}\e E^\e d\mathcal V(x)\leq  \int_{B_{R_0}(y)} \phi\frac{d^2(x)}{s-t}\e E^\e d\mathcal V(x)\\
 \leq \int_{B_{R_0}\cap \{d(x)>\sqrt[4]{s-t}\}}  \zeta\eta  \frac{d^2(x)}{s-t}\e E^\e  d\mathcal V(x) + \hspace{.7 cm} \\
 \int_{B_{R_0}\cap \{d(x)<\sqrt[4]{s-t}\}}  \zeta\eta  \frac{d^2(x)}{s-t}\e E^\e  d\mathcal
 V(x)\,.\hspace{.8 cm}
 \end{split}
\end{equation}
Furthermore, in view of \eqref{e170}, there exists a positive
constant $\underline C_{K}$, independent of $\e$, such that
\begin{equation}\label{u62b}
\mu_t^\e(K)\leq \underline C_{K}.
\end{equation}
We have
\begin{equation}\label{u63}
\int_{B_{R_0}\cap \{d(x)<\sqrt[4]{s-t}\}}  \zeta\eta
\frac{d^2(x)}{s-t}\e E^\e  d\mathcal V(x)\leq \frac{\tilde
C}{\sqrt{s-t}}\int_M \phi \e E^\e d\mathcal V(x)\,;
\end{equation}
moreover, since $B_{R_0}(y)\subset K$, from \eqref{u62b} we can
infer that
\begin{equation}\label{u64}
\begin{split}
 \int_{B_{R_0}\cap \{d(x)>\sqrt[4]{s-t}\}}  \zeta\eta  \frac{d^2(x)}{s-t}\e E^\e  d\mathcal V(x) \hspace{1.9 cm}\\
 \leq \underline C  \int_{B_{R_0}\cap \{d(x)>\sqrt[4]{s-t}\}}  \frac{d^2(x)}{(s-t)^{\frac{N+1}2}} e^{-\frac{d^2(x)}{4(s-t)}} \e E^\e
 d\mathcal V(x)\hspace{1 cm}\\
 \leq \underline C  \int_{B_{R_0}\cap \{d(x)>\sqrt[4]{s-t}\}}  \frac{d^2(x)}{s-t} e^{-\frac{d^2(x)}{8(s-t)}} \frac{e^{-\frac 1{4\sqrt{s-t}}}}{(s-t)^{\frac{N-1}2}}  \e E^\e  d\mathcal V(x)
 \leq \bar C_K\,,
\end{split}
\end{equation}
for some positive constants $\tilde C, \underline C, \bar C_K$ independent of $\e$.

On the other hand, observe that the functions $(x, t)\mapsto \eta(x,t)$ and $(x,t)\mapsto \frac{\eta(x,t)}{s-t}$
are bounded in $\left\{\frac{R_0}2 \leq d(x) \leq R_0\right\}\times (0, s)\big)$;
thus
\begin{equation}\label{u65}
\int_{\{\frac{R_0}2\leq d(x)\leq R_0\}}\e E^\e d\mathcal V(x)\leq \check C_K,
\end{equation}
for some positive constant $\check C_K$ independent of $\e$.

\smallskip

From \eqref{e32}, \eqref{u62}-\eqref{u64} we obtain \eqref{e47}, with
\[ C_3=C\tilde C_K   ,\;\; C_4=\max\{C\bar C_K, \tilde C_{R_0}\check C_K\}\,. \]

 \hfill $\square$

\medskip

From Lemma \ref{lemma5} and the asymptotic control of discrepancy
in Proposition \ref{prop2} we finally deduce the main result of this section.
\begin{theorem}\label{prop3}
Let assumption $(H_0)$ be satisfied. Let $u^\e$ be the solution to
problem \eqref{e1}-\eqref{e1a}. Suppose that \eqref{e121b} and
\eqref{e170} with $\t=0$ hold true. Let $K\subset M$ be a compact
subset, $y\in K, s>0.$ Let $\phi:=\eta\zeta$ with $\eta$ and
$\zeta$ as in Lemma \ref{lemma20}. Then for every $0<\e<1$
inequality \eqref{e60} holds true, for all $0\le t<s$,  $C_3, C_4$
being as in Lemma \ref{lemma5}, and for some positive constant
$C_5$ independent of $\e$ and $(y,s)$. As a consequence, for all
$0\leq t_0<t<s$, inequality \eqref{e61} holds true.
\end{theorem}

\noindent {\it Proof\,.\,\,}
%\begin{proof}
By \eqref{u13},
\begin{equation} \label{u66}
\int_M \frac{\eta}{s-t}\xi^\e_t d\mathcal V(x)\leq \frac{C_5}{\sqrt{s-t}},
\end{equation}
for some positive constant $C_5$, independent of $\e, y, s$. From \eqref{u66} and \eqref{e47} we can deduce \eqref{e60}.
Thus \eqref{e61} follows from
Gronwall's inequality.
%\end{proof}
\hfill $\square$

\medskip

As a consequence, the next proposition gives uniform density bounds for the measures $\mu_t^\e$ at small scales.
\begin{proposition}\label{prop12}
Let assumption $(H_0)$ be satisfied. Let $u^\e$ be the solution to
problem \eqref{e1}-\eqref{e1a}. Suppose that \eqref{e121b} and
\eqref{e170} with $\t=0$ hold true. Then, for each compact subset
$K\subset M$,
\[\int_{M}\phi_{y,s}(x,t) d\mu^\e_t(x)\le \underline C\quad \textrm{for all}\;\; y\in K,  0\le t<s\,, \leqno (G^\e_1)\]
for some $\underline C=\underline C_K>0$, and
\[\mu^\e_t(B_R(x))\le \omega_{N-1} D_0 R^{N-1}\quad \textrm{for all}\;\; x\in K, 0< R<\tilde R, t\ge 0, \leqno (G^\e_2)\]
for some $0<\tilde R<R_0$ and $D_0=D_0(\underline C, \tilde R)>0$.
\end{proposition}

\noindent{\it Proof\,.} From \eqref{e61} and \eqref{e170} it
easily follows $(G^\e_1)$. Observe that, for some $0<\tilde R<R_0,
C=C(\tilde R)>0$, there holds, for all $x,y\in M, 0<R<\tilde R,
t\geq 0$
\[\frac 1 {R^{N-1}}\chi_{B_R(y)}(x)\leq C \phi_{y,s}(x,t),\]
whenever $s-t=R^2.$ This, combined with $(G_1^\e)$, yields
$(G_2^\e)$. \hfill $\square$

%\begin{remark}\label{ossp12}
%Observe that in the proof of Lemma \ref{lemma5}, Propositions
%\ref{prop3} and \ref{prop12} we did not use the global condition
%\eqref{e162}, but we only used \eqref{e170}.
%\end{remark}

\subsection{Further compactness properties} Concerning the family
$\{u^\e\}_{0<\e<1}$ of solutions to problem \eqref{e1}-\eqref{e1a}
we have the next compactness result.
\begin{proposition}\label{prop10}
Let assumption $(H_0)$ be satisfied. Let $\{u^\e\}_{0<\e<1}$
be a family of solutions to problem \eqref{e1}-\eqref{e1a} which is uniformly bounded, i.e. satisfying \eqref{e121b}, and such that for each $T>0$ and for each compact set $K\subset M$
\begin{equation}\label{u76}
\sup_{0<\e<1} \, \sup_{t\in[0,T]} \e \int_{K} E^\e(x,t) d \mathcal{V} \leq C \, ,
\end{equation}
for some $C=C(K,T)>0$ independent of $\e$.
  Set
$\mathcal F(s):= \int_0^s \sqrt{F(\t)}d\t.$ Then there exist a
subsequence $\{\mathcal F(u^{\e_n})\}\subset \{\mathcal F(u^\e)\}$
and a function $v$ such that
\begin{equation}\label{e139}
\mathcal F(u^{\e_n})\to v\,\quad {as}\;\;n\to\infty\,,
\end{equation}
both in $C^{0,\a}_{loc}\big([0,\infty);L^1_{loc}(M)\big)$ for each
$0\le\a<\frac 1 2$ and in $BV_{loc}(M\times(0,\infty))$, and
\begin{equation}\label{e139a}
u^{\e_n}\to u^*:=\mathcal F^{-1}(v)\quad {as}\;\;n\to\infty\,,
\end{equation}
in $C^{0}_{loc}\big([0,\infty);L^1_{loc}(M)\big)$.

Moreover,
$u^*\in L^\infty_{loc}\big((0,\infty); BV_{loc}(M)\big) \cap BV_{loc}(M\times(0,\infty))$,
$$|u^*| = 1\quad \textrm{a.e. in}\;\; M\times (0,\infty),$$ and
the jump set of $u^*(\cdot, t)$ is locally $(N-1)-$ rectifiable
for a.e. $t>0$.
\end{proposition}

\noindent{\it Proof\,.} Let $T>0$ be fixed and let $\varphi \in C^\infty_0(M)$. By $(H_0)$ and \eqref{e1}-\eqref{e1a},  multiplying
\eqref{e1} by $\varphi^2\pa_t u^\e$ and integrating by parts we easily
obtain
\[\int_0^T \int_M \e \varphi^2 (\pa_t u^\e)^2 d \mathcal V dt\, \le \int_0^T \int_M \e  \varphi^2(\pa_t u^\e)^2 d \mathcal V dt\, +  \e \int_{M}\varphi^2 E^\e(x,T) d \mathcal{V}\le\]
\[\e \int_{M}\varphi^2 E^\e(x,0) d \mathcal{V}-2\e \int_0^T \int_M \varphi \pa_t u^\e \nabla u^\e \nabla \varphi d\mathcal{V} dt \,  . \]
Since $| \nabla u^\e|^2 \leq 2 E^\e(x,t)$, applying Young's inequality and \eqref{u76} with $K=supp \, \varphi$ we easily obtain for any $\varphi \in C^\infty_0(M)$ the uniform bound
\begin{equation}\label{e122}
\int_0^T \int_M \e \varphi^2 (\pa_t u^\e)^2 d \mathcal V dt\, \le C.
\end{equation}
In addition, in view of \eqref{e121b} we also have
\begin{equation}\label{e123}
|\mathcal F(u^\e)|\le C\quad \textrm{in}\;\; M\,
\end{equation}
for every $\e>0.$ From \eqref{u76}-\eqref{e123} and Young's inequality it follows that
\begin{equation}\label{u77}
\|\mathcal F(u^\e)\|_{L^{\infty}\big((0,T);BV_{loc}(M)\big)}+\| \partial_t  \mathcal{F}(u^\e) \|_{L^1 \big(0,T;L^1_{loc}(M) \big)} \le
C\;\;
\end{equation}
for every $\e>0$. As a consequence, $\mathcal{F}(u^\e)$ is $*-$weakly compact in $BV_{loc}(M\times (0,\infty))$, since $T>0$ is arbitrary.
%Furthermore, from \eqref{e122}-\eqref{e123} we
%can infer that
%\begin{equation}\label{e122b}
%\int_0^T \int_M \e (\pa_t \mathcal F(u^\e))^2 d \mathcal V dt\,
%\le C.
%\end{equation}
%for every $\e>0$.

In view of \eqref{u76},  \eqref{e123} and the estimate on the time derivative in \eqref{u77} we
can also infer that
\[\|\mathcal F(u^\e)\|_{C^{0, 1/2}\big([0,T];L^1_{loc}(M)\big)}\le
C(T)\,,\] for every $\e>0.$ So, by Ascoli-Arzelà theorem,
$\{\mathcal F(u^\e)\}_{\e>0}$ is compact in
$C^{0,\a}_{loc}\big([0,\infty);L^1_{loc}(M)\big)$ for each $0\le
\a<\frac 12$ and the same clearly holds for $\{u^\e\}_{\e>0}$ by continuity of $\mathcal{F}^{-1}$. Thus, by a diagonal argument there exist a subsequence $\{\mathcal
F(u^{\e_n})\}\subset \{\mathcal F(u^\e)\}$ and a function $v$ such
that \eqref{e139}-\eqref{e139a} holds. Moreover, it is direct to
see from the uniform bound on the potential energy in \eqref{u76} that $|u^*|=1$ a.e. in $M\times (0,\infty)$,
 and in turn that $u^*\in L^\infty\big((0,\infty);
BV_{loc}(M)\big)\cap BV_{loc}(M\times (0,\infty))$. Then the last statement follows
from rectifiability of jump set for $BV$ functions .\hfill
$\square$

\end{document}